\documentclass[12pt]{article}
\usepackage{latexsym, amssymb, amsmath, amscd, amsfonts, epsfig, graphicx, colordvi,verbatim,ifpdf,extarrows}
\usepackage{amsfonts, amsmath, amssymb}
\usepackage{amssymb,amsfonts,amsmath,latexsym,epsfig,cite, psfrag,eepic,color}
\usepackage{amscd,graphics}
\usepackage{latexsym, amssymb,  amsmath,amscd, amsfonts, epsfig, graphicx, colordvi,amsthm}

\usepackage{graphicx}
\usepackage{epstopdf}
\usepackage{color}
\usepackage{ifpdf}
\usepackage{fancybox}
\usepackage[font=small,labelfont=bf,labelsep=none]{caption}
\usepackage{float}

\newtheorem{thm}{Theorem}[section]
\newtheorem{prop}[thm]{Proposition}
\newtheorem{conj}[thm]{Conjecture}

\newtheorem{defi}[thm]{Definition}
\newtheorem{lem}[thm]{Lemma}

\def\pf{\noindent{\it Proof.} }
\def\qed{\nopagebreak\hfill{\rule{4pt}{7pt}}
\medbreak}

\numberwithin{equation}{section}

\def\qed{\nopagebreak\hfill{\rule{4pt}{7pt}}
\medbreak}

\newlength{\boxedparwidth}
  {\begin{center} \begin{tabular}{|@{\hspace{.315in}}c@{\hspace{.15in}}|}
                  \hline \\ \begin{minipage}[t]{\boxedparwidth}
                  \setlength{\parindent}{.25in}}%
  {\end{minipage} \\ \\ \hline \end{tabular} \end{center}}

\begin{document}

\begin{center}

 {\Large \bf Overpartitions and Bressoud's conjecture, II}
\end{center}

\begin{center}
 {Thomas Y. He}$^{1}$, {Kathy Q. Ji}$^{2}$  and
  {Alice X.H. Zhao}$^{3}$ \vskip 2mm

$^{1}$ School of Mathematical Sciences,  Sichuan Normal University, Sichuan 610068, P.R. China\\[6pt]

$^{2}$ Center for Applied Mathematics,  Tianjin University, Tianjin 300072, P.R. China\\[6pt]

  $^{3}$ College of Science, Tianjin University of Technology, Tianjin 300384, P.R. China

   \vskip 2mm

    $^1$heyao@sicnu.edu.cn, $^2$kathyji@tju.edu.cn,  $^3$zhaoxh@email.tjut.edu.cn
\end{center}

\vskip 6mm   {\bf Abstract.}  The main objective of this paper is to present an answer to Bressoud's conjecture for the case $j=0$,  resulting in a complete solution to  Bressoud's conjecture. The case for $j=1$ has been recently resolved by Kim. Using the connection  established in our previous paper between the ordinary partition function $B_0$ and the overpartition function $\overline{B}_1$, we found that the proof of  Bressoud's conjecture for the case $j=0$ is equivalent to establishing an overpartition analogue of the conjecture for the case $j=1$. By generalizing Kim's method, we obtain the desired  overpartition analogue of Bressoud's conjecture for the case $j=1$, which eventually enables us to confirm Bressoud's conjecture for the case $j=0$.



\section{Introduction}

This is the second in a series of papers addressing Bressoud's conjecture. In 1980,  Bressoud \cite{Bressoud-1980} put forward a conjecture for a   general partition identity that implies many classical results in
the theory of partitions, such as Euler's partition theorem,  the Rogers-Ramanujan-Gordon identities, the Andrews-G\"ollnitz-Gordon identities and so on.  To state Bressoud's conjecture, let us  recall some common notation and terminology on partitions from \cite[Chapter 1]{Andrews-1976}. A partition $\pi$ of a positive integer $n$ is a finite non-increasing sequence of positive integers $\pi=(\pi_1,\pi_2,\ldots,\pi_\ell)$ such that $\sum_{i=1}^{\ell} \pi_i=n$.  The weight of $\pi$ is the sum of its parts, denoted $|\pi|$.

 Throughout this paper, we assume that $\alpha_1,\alpha_2,\ldots, \alpha_\lambda$ and $\eta$ are integers such that
\begin{equation*}\label{cond-alpha}
0<\alpha_1<\alpha_2<\cdots<\alpha_\lambda<\eta, \quad \text{and} \quad \alpha_i=\eta-\alpha_{\lambda+1-i}\quad \text{for} \quad 1\leq i\leq \lambda.
\end{equation*}

 Bressoud  \cite{Bressoud-1980} introduced the following two partition functions.

\begin{defi}[Bressoud] For $j=0$ or $1$  and $(2k+j)/2> r\geq\lambda\geq0$,   define the partition function $A_j(\alpha_1,\ldots,\alpha_\lambda;\eta,k,r;n)$ to be the number of partitions of $n$ into parts congruent to $0,\alpha_1,\ldots,\alpha_\lambda\pmod\eta$ such that
\begin{itemize}
\item[\rm{(1)}] If $\lambda$ is even, then only multiples of $\eta$ may be repeated and no part is congruent to $0,\pm\eta(r-\lambda/2)  \pmod{\eta(2k-\lambda+j)}${\rm{;}}

    \item[\rm{(2)}] If $\lambda$ is odd and $j=1$, then only multiples of ${\eta}/{2}$ may be repeated, no part is congruent to $\eta\pmod{2\eta}$, and no part is congruent to $0,\pm{\eta}(2r-\lambda)/{2} \pmod {\eta(2k-\lambda+1)}${\rm{;}}

        \item[\rm{(3)}] If $\lambda$ is odd and $j=0$, then only multiples of ${\eta}/{2}$ which are not congruent to ${\eta}(2k-\lambda)/{2}\pmod{\eta(2k-\lambda)}$ may be repeated, no part is congruent to $\eta\pmod{2\eta}$, no part is congruent to $0\pmod{2\eta(2k-\lambda)}$, and no part is congruent to $\pm{\eta}(2r-\lambda)/{2} \pmod {\eta(2k-\lambda)}$.
  \end{itemize}
  \end{defi}

\begin{defi}[Bressoud] \label{Bress-B-function} For $j=0$ or $1$  and $k\geq r\geq\lambda\geq0$,  define the partition function
$B_j(\alpha_1,\ldots,\alpha_\lambda;\eta,k,r;n)$  to be the number of partitions  $\pi=(\pi_1,\pi_2,\ldots,\pi_\ell)$  of $n$  satisfying the following conditions{\rm{:}}
\begin{itemize}
\item[{\rm (1)}] For $1\leq i\leq \ell$, $\pi_i\equiv0,\alpha_1,\ldots,\alpha_\lambda\pmod{\eta}${\rm{;}}

\item[{\rm (2)}] Only multiples of $\eta$ may be repeated{\rm{;}}

\item[{\rm (3)}] For $1\leq i\leq \ell-k+1$, $ \pi_i\geq\pi_{i+k-1}+\eta$  with strict inequality if $\eta\mid\pi_i${\rm{;}}

\item[{\rm (4)}] At most $r-1$ of the $\pi_i$ are less than or equal to $\eta${\rm{;}}

\item[{\rm (5)}] For $1\leq i\leq \ell-k+2$, if $\pi_i\leq\pi_{i+k-2}+\eta$ with strict inequality if $\eta\nmid\pi_i$,  then \[[\pi_i/\eta]+\cdots+[\pi_{i+k-2}/\eta]\equiv r-1+V_\pi(\pi_i)\pmod{2-j},\]
  where $V_\pi(N)$ denotes the number of parts not exceeding $N$ which are not divisible by $\eta$ in $\pi$ and $[\ ]$ denotes the greatest integer function.
  \end{itemize}
 \end{defi}

 Bressoud's conjecture can be stated as follows.

  \begin{conj}[Bressoud] \label{Bressoud-conjecture-j} For $j=0$ or $1$, $(2k+j)/2> r\geq\lambda\geq0$ and $n\geq 0$,
 \begin{equation*}\label{Bressoud-conj-1}
A_j(\alpha_1,\ldots,\alpha_\lambda;\eta,k,r;n)=B_j(\alpha_1,\ldots,\alpha_\lambda;\eta,k,r;n).
 \end{equation*}
\end{conj}

Bressoud's conjecture was known in some special cases, see,  Andrews \cite{Andrews-1974m}, Bressoud  \cite{Bressoud-1980} and Kim and Yee \cite{Kim-Yee-2014}.   The general case for $j=1$ was recently resolved by Kim  \cite{Kim-2018}. The main objective of this paper is to present an answer to Bressoud's conjecture for the case $j=0$, resulting in a complete solution to  the conjecture.
It turns out that the overpartition analogues of the partition
 functions ${A}_1(\alpha_1,\ldots,\alpha_\lambda;\eta,
 k,r;n)$ and  ${B}_1(\alpha_1,\ldots,\alpha_\lambda;\eta,
 k,r;n)$ introduced in our previous paper \cite{he-ji-zhao}
 play an important role in the proof of Conjecture \ref{Bressoud-conjecture-j} for the case $j=0$.
An overpartition, introduced by Corteel and Lovejoy \cite{Corteel-Lovejoy-2004},  is a partition such that the first occurrence of a part can be overlined. In fact, they arise in the contexts of combinatorics \cite{Bessenrodt-Pak-2004, Pak-2006}, $q$-series \cite{Corteel-Hitczenko-2004}, the theory of symmetric functions \cite{Brenti-1999, Desrosiers-Lapointe-Mathieu-2006}, representation theory \cite{Kang-Kwon-2004} and mathematical physics \cite{Fortin-Jacob-Mathieu-2005, Fortin-Jacob-Mathieu-2005b}. They are also known by different names such as standard MacMahon diagrams, joint partitions, jagged partitions and dotted partitions.

In \cite{he-ji-zhao}, we introduced the following two functions $\overline{A}_1(\alpha_1,\ldots,\alpha_\lambda;\eta,k,r;n)$ and $\overline{B}_1(\alpha_1,\break$$\ldots,\alpha_\lambda;\eta,k,r;n)$ defined on the set of overpartitions.

\begin{defi}{\rm \!\!\!\cite[Definition 1.15]{he-ji-zhao}}	
For  $k> r\geq\lambda\geq0$,  define the partition function $\overline{A}_1(\alpha_1,\ldots,\alpha_\lambda;\eta,k,r;n)$ to be the number of overpartitions of $n$ into parts congruent to $0,\alpha_1,\ldots,\alpha_\lambda\pmod\eta$ such that
	\begin{itemize}
		\item[\rm{(1)}] If $\lambda$ is even, then only multiples of $\eta$ may be non-overlined and  there is no non-overlined part congruent to $0,\pm\eta(r-\lambda/2) \pmod {\eta(2k-\lambda)}${\rm{;}}
		
		\item[\rm{(2)}] If $\lambda$ is odd, then only multiples of ${\eta}/{2}$ may be non-overlined, no non-overlined part is congruent to ${\eta}(2k-\lambda)/{2}\pmod{\eta(2k-\lambda)}$, no non-overlined part is congruent to $\eta \pmod{2\eta}$, no non-overlined part is congruent to $0 \pmod{2\eta(2k-\lambda)}$, no non-overlined part is congruent to $\pm{\eta}(2r-\lambda)/{2} \pmod {\eta(2k-\lambda)}$, and no overlined part is congruent to ${\eta}/{2}\pmod \eta$ and not congruent to ${\eta}(2k-\lambda)/{2}\pmod{\eta(2k-\lambda)}${\rm{.}}

	\end{itemize}
\end{defi}

 \begin{defi}{\rm\!\!\! \cite[Definition 1.14]{he-ji-zhao}} \label{defi-O-B} For $k\geq r\geq \lambda\geq0$, define the partition function $\overline{B}_1(\alpha_1,\ldots,\alpha_\lambda;\eta,k,r;n)$ to be the number of overpartitions $\pi=(\pi_1,\pi_2,\ldots,\pi_\ell)$ of $n$  subject to  the following conditions{\rm{:}}
 \begin{itemize}
  \item[{\rm (1)}] For $1\leq i\leq \ell$, $\pi_i\equiv0,\alpha_1,\ldots,\alpha_\lambda\pmod{\eta}${\rm{;}}

 \item[{\rm (2)}] Only  multiples of $\eta$ may be non-overlined{\rm{;}}

 \item[{\rm (3)}] For $1\leq i\leq \ell-k+1$, $\pi_i\geq\pi_{i+k-1}+\eta$ with strict inequality if  $\pi_i$  is non-overlined{\rm{;}}

 \item[{\rm (4)}] At most $r-1$ of the $\pi_i$ are less than or equal to $\eta${\rm{.}}
\end{itemize}
\end{defi}
 Here and in the sequel, we adopt the following convention:    For a positive integer $t$,  we define $t\pm \eta$ (resp. $\overline{t}\pm \eta$) as a non-overlined part of size  $t\pm \eta$ (resp. an overlined part of size ${t\pm \eta}$). We impose the following order on the parts of an overpartition:
 \[1<\bar{1}<2<\bar{2}<\cdots.\]
As   mentioned in \cite{he-ji-zhao}, we say that $\overline{B}_1(\alpha_1,\ldots,\alpha_\lambda;\eta,k,r;n)$ (resp. $\overline{A}_{1}(\alpha_1,\ldots,\alpha_\lambda;\eta,k,r;n)$) can be considered as an overpartition analogue of  ${B}_1(\alpha_1,\ldots,\alpha_\lambda;\eta,k,r;n)$ (resp.  ${A}_0(\alpha_1,\ldots,\break$$\alpha_\lambda;\eta,k,r;n)$) because    for an overpartition $\pi$ counted by $\overline{B}_1(\alpha_1,\ldots,\alpha_\lambda;\eta,k,r;n)$ (resp. $\overline{A}_{1}(\alpha_1,\ldots,\alpha_\lambda;\eta,k,r;n)$) without overlined parts divisible by $\eta$, if we change the overlined parts in $\pi$ to non-overlined parts, then we get an ordinary partition counted by ${B}_1(\alpha_1,\ldots,\alpha_\lambda;\eta,k,r;n)$ (resp. ${A}_{0}(\alpha_1,\ldots,\alpha_\lambda;\eta,k,r;n)$).

One of the main results in our previous paper \cite{he-ji-zhao} is   the following connection between $\overline{B}_1(\alpha_1,\ldots,\alpha_\lambda;
	\eta,k,r;n)$ and   ${B}_0(\alpha_1,\ldots,\alpha_\lambda;
	\eta,k,r;n)$, which is useful in the proof of Conjecture \ref{Bressoud-conjecture-j} for the case $j=0$.

\begin{thm}{\rm\!\!\! \cite[Theorem 1.16]{he-ji-zhao}} \label{rel-over1}For   $k\geq r\geq \lambda\geq0$ and $k>\lambda$,
	\begin{equation*}\label{new-b-0-over}
	\begin{split}
	&\sum_{n\geq0}\overline{B}_1(\alpha_1,\ldots,\alpha_\lambda;
	\eta,k,r;n)q^n=(-q^\eta;q^\eta)_\infty
	\sum_{n\geq0}B_0(\alpha_1,\ldots,\alpha_\lambda;
	\eta,k,r;n)q^n.
	\end{split}
	\end{equation*}
\end{thm}

The generating function of $\overline{A}_1(\alpha_1,\ldots,\alpha_\lambda;\eta,k,r;n)$ is also established in our previous paper \cite{he-ji-zhao}.

 \begin{thm}{\rm \!\!\! \cite[Theorem 1.18]{he-ji-zhao}}\label{gf-overline-a}
 For  $k>r\geq\lambda\geq0$,
\begin{equation*}\label{overpartition-Afunction-e}
\begin{split}
&\sum_{n\geq0}\overline{A}_1(\alpha_1,\ldots,\alpha_\lambda;\eta,k,r;n)q^n\\[5pt]
&\quad=\frac{(-q^{\alpha_1},\ldots,-q^{\alpha_\lambda},-q^{\eta};q^{\eta})_\infty(q^{\eta(r-\frac{\lambda}{2})},
	q^{\eta(2k-r-\frac{\lambda}{2})},
	q^{\eta(2k-\lambda)};
	q^{\eta(2k-\lambda)})_\infty}{(q^\eta;q^\eta)_\infty}.
\end{split}
\end{equation*}
\end{thm}

Throughout this paper, we assume that $|q|<1$ and employ the standard notation:
\[(a;q)_\infty=\prod_{i=0}^{\infty}(1-aq^i), \quad (a;q)_n=\frac{(a;q)_\infty}{(aq^n;q)_\infty},\]
and
\[(a_1,a_2,\ldots,a_m;q)_\infty=(a_1;q)_\infty(a_2;q)_\infty\cdots(a_m;q)_\infty.
\]

It should be noted that Bressoud \cite{Bressoud-1980} obtained the following generating function of $A_0(\alpha_1,\ldots,\alpha_\lambda;\eta,k,r;n)$.
\begin{thm}[Bressoud] \label{thm-gf-a}For $k> r\geq\lambda\geq0$,
 \begin{equation*}\label{Bressoud-conj-defi-e}
 \begin{split} &\sum_{n\geq0}A_0(\alpha_1,\ldots,\alpha_\lambda;\eta,k,r;n)q^n\\[5pt]
 &\quad=\frac{(-q^{\alpha_1},\ldots,-q^{\alpha_\lambda};q^{\eta})_\infty
 (q^{\eta(r-\frac{\lambda}{2})},q^{\eta(2k-r-\frac{\lambda}{2})},
 q^{\eta(2k-\lambda)};q^{\eta(2k-\lambda)})_\infty}{(q^\eta;q^\eta)_\infty}.
 \end{split}
 \end{equation*}
\end{thm}

In this paper, we establish the following generating function of $\overline{B}_1(\alpha_1,\ldots,\alpha_\lambda;\eta,k,r;n)$ by generalizing Kim's  method in \cite{Kim-2018}.

 \begin{thm} \label{Over-bre-1a} For $k> r\geq \lambda\geq0$,
\begin{equation}\label{proof-0-A}
 \begin{split}
 &\sum_{n\geq0}\overline{B}_1(\alpha_1,\ldots,\alpha_\lambda;\eta,k,r;n)q^n\\
 &\quad=
\frac{{(-q^{\alpha_1},\ldots,-q^{\alpha_\lambda},-q^\eta;q^\eta)_\infty}(q^{\eta(r-\frac{\lambda}{2})},
 q^{\eta(2k-r-\frac{\lambda}{2})},
 q^{\eta(2k-\lambda)};
 q^{\eta(2k-\lambda)})_\infty}{(q^\eta;q^\eta)_\infty}.
 \end{split}
 \end{equation}
 \end{thm}

Applying Theorem \ref{Over-bre-1a} to Theorem \ref{rel-over1}, and according to  Theorem \ref{thm-gf-a}, we conclude that  Conjecture \ref{Bressoud-conjecture-j}  holds for the case $j=0$.

 Combining Theorem \ref{gf-overline-a} and Theorem \ref{Over-bre-1a}, we also get the following partition identity, which can be viewed as  an overpartition analogue of Bressoud's conjecture for $j=1$.

\begin{thm}\label{Over-bre-1}
For  $k> r\geq \lambda\geq0$ and $n\geq 0,$
  \begin{equation*}\label{main}
 \overline{A}_1(\alpha_1,\ldots,\alpha_\lambda;\eta,k,r;n)=\overline{B}_1(\alpha_1,\ldots,\alpha_\lambda;\eta,k,r;n).
 \end{equation*}
\end{thm}

Through the utilization of Theorem \ref{Over-bre-1a} and the application of Bailey pairs, we can formulate the generating function of $\overline{B}_1(\alpha_1,\ldots,\alpha_\lambda;\eta,k,r;n)$ as the following multi-summation identity.

\begin{thm} \label{Over-bre-1bbb} For $k> r> \lambda\geq0,$
\begin{equation*}
 \begin{split}
 &\sum_{n\geq0}\overline{B}_1(\alpha_1,\ldots,\alpha_\lambda;\eta,k,r;n)q^n\\[5pt]
 &\quad=\sum_{N_1\geq\cdots\geq N_{k-1}\geq0}
 \frac{q^{\eta(N_{1}^{2}+\cdots+N_{k-1}^{2}+ N_r+\cdots+N_{k-1})}(1+q^{-\eta N_{r}})(-q^{\eta-\eta N_{\lambda+1}};q^{\eta})_{N_{\lambda+1}-1}
 (-q^{\eta+\eta N_{\lambda}};q^{\eta})_{\infty}}
 {(q^{\eta};q^{\eta})_{N_{1}-N_{2}}\cdots(q^{\eta};q^{\eta})_{N_{k-2}-N_{k-1}}(q^{\eta}
 ;q^{\eta})_{N_{k-1}}}\\[5pt]
 &\qquad\qquad\times
 \prod_{s=1}^\lambda(-q^{\eta-\alpha_{s}-\eta N_{s}};q^{\eta})_{N_{s}}
\prod_{s=2}^\lambda(-q^{\eta-\alpha_{s}+\eta N_{s-1}};q^{\eta})_{\infty}.
 \end{split}
 \end{equation*}
 \end{thm}
 It would be interesting to give a combinatorial proof of Theorem \ref{Over-bre-1bbb}.

This article is organized as follows. In Section 2, we aim to prove Theorem \ref{G-B-O-1} with the aid of Bailey pairs. As a result, Theorem \ref{Over-bre-1bbb} follows directly from the conjunction of Theorem \ref{Over-bre-1a} and Theorem \ref{G-B-O-1}.  Section 3 is dedicated to proving Theorem \ref{Over-bre-1a}. Initially, we establish that proving Theorem \ref{Over-bre-1a} is sufficient to show  the validity of Theorem \ref{lambda-r}. Following this, we present an outline proof of Theorem \ref{lambda-r}, which equivalently demonstrates the combinatorial statement in Theorem \ref{lambdathm}.  In order to establish the desired bijection as stated in Theorem \ref{lambdathm}, we begin by revisiting the definitions of the Gordon marking and the reverse Gordon marking of an overpartition counted by $\overline{{B}}_1(\alpha_1,\ldots,\alpha_\lambda;\eta,k,r;n)$, as previously established in our prior work \cite{he-ji-zhao}. Subsequently, we review the $(k-1)$-addition and its inverse map (i.e., the $(k-1)$-subtraction) introduced by Kim \cite{Kim-2018} in the context of overpartitions. Additionally, we recall  the $(k-1)$-insertion and its inverse map (i.e., the $(k-1)$-separation) defined in our prior work \cite{he-ji-zhao}.  These operations allow us to provide the desired bijection  presented in Theorem \ref{lambdathm}.  In Section 4, we provide an example for the illustration of the bijection in the proof of Theorem \ref{lambdathm}.

\section{Proof of Theorem \ref{Over-bre-1bbb}}

The main objective of this section is to give a proof of Theorem \ref{G-B-O-1} by using  Bailey pairs. Consequently, Theorem \ref{Over-bre-1bbb} is immediately derived from the combination  of Theorem \ref{Over-bre-1a} and Theorem \ref{G-B-O-1}.

\begin{thm}\label{G-B-O-1}For $k\geq r>\lambda \geq 0$, we have
\begin{equation*}\label{G-B-O-1-eq}
  \begin{split}
  &\sum_{N_1\geq\cdots\geq N_{k-1}\geq0}
 \frac{q^{\eta(N_{1}^{2}+\cdots+N_{k-1}^{2}+ N_r+\cdots+N_{k-1})}(1+q^{-\eta N_{r}})(-q^{\eta-\eta N_{\lambda+1}};q^{\eta})_{N_{\lambda+1}-1}
 (-q^{\eta+\eta N_{\lambda}};q^{\eta})_{\infty}}
 {(q^{\eta};q^{\eta})_{N_{1}-N_{2}}\cdots(q^{\eta};q^{\eta})_{N_{k-2}-N_{k-1}}(q^{\eta}
 ;q^{\eta})_{N_{k-1}}}\\[5pt]
 &\qquad\times
 \prod_{s=1}^\lambda(-q^{\eta-\alpha_{s}-\eta N_{s}};q^{\eta})_{N_{s}}
\prod_{s=2}^\lambda(-q^{\eta-\alpha_{s}+\eta N_{s-1}};q^{\eta})_{\infty}\\[5pt]
&\qquad\qquad=\frac{(-q^{\alpha_{1}},\ldots,-q^{\alpha_{\lambda}},-q^{\eta};q^{\eta})_{\infty}(q^{(r-\frac{\lambda}{2})\eta},
q^{(2k-r-\frac{\lambda}{2})\eta}
  ,q^{(2k-\lambda)\eta};q^{(2k-\lambda)\eta})_{\infty}}
{(q^{\eta};q^{\eta})_{\infty}},
  \end{split}
\end{equation*}
where we assume that $N_{k}=0$.
\end{thm}
 It should be noted that the proof of Theorem \ref{G-B-O-1} is much similar to the proof of Theorem 1.8 in \cite{He-Ji-Wang-Zhao}. For more information on Bailey pairs, see, for example, \cite{Agarwal-Andrews-Bressoud, Andrews-1986, Andrews-2000, Bressoud-Ismail-Stanton-2000, Lovejoy-2004b, Paule-1987,Warnaar-2001}.
Recall that  a pair of sequences $(\alpha_n(a,q),\beta_n(a,q))$ is called a Bailey pair relative to $(a,q)$ (or  a Bailey pair for short) if  for  $n\geq 0,$
\begin{equation*}\label{bailey pair}
\beta_n(a,q)=\sum_{r=0}^n\frac{\alpha_r(a,q)}{(q;q)_{n-r}(aq;q)_{n+r}}.
\end{equation*}
When $k>r+1\geq 2$, it turns out that  the proof of Theorem \ref{G-B-O-1} reduces to applying the Bailey pairs stated in Lemma \ref{BPG} to the relation stated in Proposition \ref{lc}.

\begin{lem} \label{BPG}  {\rm  \cite[(2.9)]{He-Ji-Wang-Zhao}}  For $k>r+1\geq 2$,
\begin{equation}{\label{BPG-eq}}
\begin{split}
\alpha_n(1,q)&=\left\{
             \begin{array}{ll}
             1, & \hbox{if $n=0$}, \\[5pt]
             (-1)^nq^{\frac{2k-2r+1}{2}n^2}(q^{\frac{2k-2r-1}{2} n}+q^{-\frac{2k-2r+1}{2} n})
             (1+q^{n})/2, & \hbox{if $n\geq1,$}
             \end{array}
               \right.\\[10 pt]
\beta_{n}(1,q)&=\sum_{n\geq N_{r+1}\geq \cdots\geq N_{k-1}\geq 0}\frac{(1+q^{n})q^{N_{r+1}^2+\cdots+{N_{k-1}^2}
+N_{r+1}+\cdots+ {N_{k-1}}}}{2(q;q)_{n-N_{r+1}}(q;q)_{N_{r+1}-N_{r+2}}\cdots
{(q;q)_{N_{k-1}}}}
\end{split}
\end{equation}
is a Bailey pair relative to $(1,q)$.

\end{lem}

\begin{prop}\label{lc}{\rm \cite[Proposition 6.6]{he-ji-zhao}} If $(\alpha_n(1,q^\eta),\beta_n(1,q^\eta))$ is a Bailey pair relative to $(1,q^\eta)$, then for $r>\lambda\geq0$,
\begin{align}\nonumber
& \sum_{n=0}^\infty\frac{2q^{(r-\frac{\lambda +1}{2})\eta n^2+\frac{\lambda +1}{2}\eta n
-(\alpha_1+\cdots+\alpha_{\lambda})n}(-q^{\alpha_1};q^{\eta})_n\cdots(-q^{
\alpha_{\lambda}};q^{\eta})_n}{(1+q^{\eta n})(-q^{\eta-\alpha_1};q^{\eta})
_n\cdots(-q^{\eta-\alpha_{\lambda}};q^{\eta})_n}\alpha_n(1,q^\eta)
\\[5pt]\nonumber
&\quad=\frac{(q^{\eta};q^{\eta})_{\infty}}{(-q^{\eta-\alpha_1};q^{\eta})_
{\infty}}\sum_{N_1\geq N_2\geq\cdots\geq N_r\geq 0}
\frac{q^{
\eta(N_{\lambda+2}^2+\cdots+N_{r}^2)+\eta\left({N_1+1 \choose2}+\cdots+{N_{\lambda+1}+1 \choose2}\right)
-(\alpha_1N_1+\cdots+\alpha_{\lambda}N_{\lambda})
}}{(q^{\eta};q^{\eta})_{N_1-N_{2}}\cdots(q^{\eta};q^{\eta})_{
N_{r-1}-N_{r}}}\\[5pt] \label{lem-cor}
&\qquad\times\frac{(-1;q^{\eta})_{N_{\lambda+1}
}(-q^{\alpha_1};q^{\eta})_{N_1}\cdots(-q^{\alpha_{\lambda
}};q^{\eta})_{N_{\lambda}}}
{(-q^{\eta};q^{\eta})_{N_{\lambda}}
(-q^{\eta-\alpha_2};q^
{\eta})_{N_1}\cdots(-q^{\eta-\alpha_{\lambda}};q^{\eta})_{
N_{\lambda-1}}}\beta_{N_r}(1,q^\eta),
\end{align}
where we assume that $N_{r+1}=0$.
\end{prop}

To prove that Theorem \ref{G-B-O-1} is valid when $k=r\geq 1$ and $k=r+1\geq 2$, we also need to apply the following proposition to generate the desired Bailey pairs.

\begin{prop}\label{bl4} {\rm \cite[Corollary 2.4]{He-Ji-Wang-Zhao}}
If $A$ is any real number and  $(\alpha_n(1,q),\beta_n(1,q))$ is a Bailey pair, where
\[\alpha_n(1,q)=\left\{
             \begin{array}{ll}
               1, & \hbox{if $n=0$}, \\[5pt]
               (-1)^nq^{An^2}(q^{(A-1)n}+q^{-(A-1)n}), & \hbox{if $n\geq 1$,}
             \end{array}
           \right.
\]
then $(\alpha_n'(1,q),\beta_n'(1,q))$ is also a Bailey pair, where
\begin{align*}
 \alpha_n'(1,q)&=\left\{
             \begin{array}{ll}
               1, & \hbox{if $n=0$}, \\[3pt]
               (-1)^nq^{An^2}(q^{(A-1)n}+q^{-An})(1+q^n)/2, & \hbox{if $n\geq 1$,}
             \end{array}
           \right. \\[10pt]
 \beta_n'(1,q)&=\beta_n(1,q)(1+q^n)/2.
\end{align*}
\end{prop}
With these consequences in hand, we are ready to prove Theorem \ref{G-B-O-1}.

\noindent{\it Proof of Theorem \ref{G-B-O-1}}: We consider the following two cases:

Case 1. When $k>r>\lambda\geq 0$. We first show that the pair of sequences $(\alpha_n(1,q),\beta_n(1,q))$ stated in \eqref{BPG-eq} is also a Bailey pair when $k=r+1\geq 2$. In this case,  we assume that $N_k=0$. We start with the following Bailey pair appearing in  \cite[B(1)]{Slater-1952},
\begin{align*}
& \alpha_n^{(0)}(1, q)= \begin{cases}1, & \text { if } n=0, \\
(-1)^n q^{3n^2 / 2}\left(q^{-n / 2}+q^{n / 2}\right), & \text { if } n \geq 1,\end{cases} \\[5pt]
& \beta_n^{(0)}(1, q)= \frac{1}{(q;q)_n}.
\end{align*}
Applying Proposition 2.4  with $A={3}/{2}$ to  ($\alpha_n^{(0)}(1, q),\,\beta_n^{(0)}(1, q)$)  yields
\begin{align*}
& \alpha_n(1, q)= \begin{cases}1, & \text { if } n=0, \\
(-1)^n q^{3n^2 / 2}\left(q^{n / 2}+q^{-3n / 2}\right)(1+q^n)/2, & \text { if } n \geq 1,\end{cases} \\[5pt]
& \beta_n(1, q)= \frac{1+q^n}{2(q;q)_n}.
\end{align*}
It follows that Lemma \ref{BPG} also holds when $k=r+1\geq 2$.

Substitute the Bailey pair \eqref{BPG-eq} from Lemma \ref{BPG}, including the case $k=r+1\geq 2$ with $q$ replaced by $q^\eta$, into \eqref{lem-cor}. Given the assumption $\alpha_i+\alpha_{\lambda+1-i}=\eta$ for $1\leq i\leq \lambda$, the simplification of the left-hand side of \eqref{lem-cor} yields
\begin{align}\nonumber
   & 1+\sum_{n=1}^\infty\frac{(-q^{\alpha_1};q^{\eta})_n\cdots(-q^{
\alpha_{\lambda}};q^{\eta})_n}{(-q^{\eta-\alpha_1};q^{\eta})
_n\cdots(-q^{\eta-\alpha_{\lambda}};q^{\eta})_n}\\[5pt] \nonumber
    & \quad\times (-1)^nq^{(k-\frac{\lambda}{2})\eta n^2+\frac{\lambda+1}{2}\eta n
-(\alpha_1+\cdots+\alpha_{\lambda})n}(q^{\frac{2k-2r-1}{2}\eta n}+q^{-\frac{2k-2r+1}{2}\eta n})\\ \nonumber
&\qquad=1+\sum_{n=1}^{\infty}(-1)^nq^{(k-\frac{\lambda}{2})\eta n^2}(q^{(k-r)\eta n}+q^{-(k-r)\eta n})\\ \label{eq-L}
&\quad\qquad=(q^{(r-\frac{\lambda}{2})\eta},q^{(2k-r-
\frac{\lambda}{2})\eta},q^{(2k-\lambda)\eta};q^{(2k-\lambda)\eta})_{\infty},
\end{align}
where the last equality follows from Jacobi's triple product identity \cite[Theorem 2.8]{Andrews-1976}.

In this case,  the  right-hand side of \eqref{lem-cor} becomes
\begin{align} \nonumber
&\frac{(q^{\eta};q^{\eta})_{\infty}}{(-q^{\eta-\alpha_1};q^{\eta})_
{\infty}}\sum_{N_1\geq \cdots\geq N_{k-1}\geq0}\frac{q^{\eta(N_{\lambda+2}^2+\cdots+N_{k-1}^2+N_r+\cdots+N_{k-1})}(1+q^{-\eta N_r})}
{(q^{\eta};q^{\eta})_{N_1-N_2}\cdots(q^{\eta};q^{\eta})_{N_{
k-2}-N_{k-1}}(q^{\eta};q^{\eta})_{N_{k-1}}}\\[5pt] \label{eq-R1}
&\hskip1cm\times\frac{q^{\eta\left({N_1+1 \choose2}+\cdots+{N_{\lambda+1}+1 \choose2}\right)-(\alpha_1N_1+\cdots+\alpha_{\lambda}N_{\lambda})}
(-q^{\eta};q^{\eta})_{N_{\lambda+1}
-1}(-q^{\alpha_1};q^{\eta})_{N_1}\cdots(-q^{\alpha_{\lambda
}};q^{\eta})_{N_{\lambda}}}{(-q^{\eta};q^{\eta})_{
N_{\lambda}}(-q^{\eta-\alpha_2};q^{\eta})_{N_1}
\cdots(-q^{\eta-\alpha
_{\lambda}};q^{\eta})_{N_{\lambda-1}}},
\end{align}
where for $N_{\lambda+1}\geq 0$,
\begin{equation}\label{sim-3} \frac{(-1;q^\eta)_{N_{\lambda+1}}}{2}=(-q^\eta;q^\eta)_{N_{\lambda+1}-1}.
\end{equation}
By utilizing the following two relations
\begin{equation}\label{sim-1}
(-q^{r};q^\eta)_{n}=q^{r n+\eta {n\choose 2}}(-q^{\eta-r-n\eta};q^\eta)_{n}
\end{equation}
and
\begin{equation}\label{sim-2}
\frac{1}{(-q^{\eta-r};q^\eta)_n}=\frac{(-q^{\eta-r+n\eta};q^\eta)_\infty}
{(-q^{\eta-r};q^\eta)_\infty},
\end{equation}
we derive that   \eqref{eq-R1} can be transformed into
\begin{align}\nonumber
&\frac{(q^{\eta};q^{\eta})_{\infty}}{(-q^{\eta-\alpha_1};q^{\eta})_
{\infty}}\sum_{N_1\geq\cdots\geq N_{k-1}\geq0}\frac{q^{\eta(N_1^2+\cdots+N_{k-1}^2+N_r+\cdots+N_{k-1})}
(1+q^{-\eta N_r})(-q^{\eta-\eta N_{\lambda+1}};q^{\eta})_{N_{\lambda+1}-1}}
{(q^{\eta};q^{\eta})_{N_1-N_2}\cdots(q^{\eta};q^{\eta})_
{N_{k-2}-N_{k-1}}(q^{\eta};q^{\eta})_{N_{k-1}}}\\[5pt] \label{eq-R2}
&\hskip1cm\times\frac{(-q^{\eta+\eta N_{\lambda}};q^{\eta})_{\infty}\prod_{s=1}^{\lambda}(-q^{\eta-\alpha_s
-\eta N_{s}};q^{\eta})_{N_s}\prod_{s=2}^{\lambda}(-q^{\eta-\alpha_s+
\eta N_{s-1}};q^{\eta})_{\infty}}{(-q^\eta;q^\eta)_\infty\prod_{s=2}^{\lambda}(-q^{
\eta-\alpha_s};q^{\eta})_{\infty}}.
\end{align}

Combining \eqref{eq-L}  and \eqref{eq-R2}, we have
\begin{align}\nonumber
&\frac{(q^{\eta};q^{\eta})_{\infty}}{(-q^{\eta-\alpha_1};q^{\eta})_
{\infty}}\sum_{N_1\geq\cdots\geq N_{k-1}\geq0}\frac{q^{\eta(N_1^2+\cdots+N_{k-1}^2+N_r+\cdots+N_{k-1})}
(1+q^{-\eta N_r})(-q^{\eta-\eta N_{\lambda+1}};q^{\eta})_{N_{\lambda+1}-1}}
{(q^{\eta};q^{\eta})_{N_1-N_2}\cdots(q^{\eta};q^{\eta})_
{N_{k-2}-N_{k-1}}(q^{\eta};q^{\eta})_{N_{k-1}}}\\[5pt] \nonumber
&\quad\times\frac{(-q^{\eta+\eta N_{\lambda}};q^{\eta})_{\infty}\prod_{s=1}^{\lambda}(-q^{\eta-\alpha_s
-\eta N_{s}};q^{\eta})_{N_s}\prod_{s=2}^{\lambda}(-q^{\eta-\alpha_s+
\eta N_{s-1}};q^{\eta})_{\infty}}{(-q^\eta;q^\eta)_\infty\prod_{s=2}^{\lambda}(-q^{
\eta-\alpha_s};q^{\eta})_{\infty}}\\ \nonumber
&\qquad=(q^{(r-\frac{\lambda}{2})\eta},q^{(2k-r-
\frac{\lambda}{2})\eta},q^{(2k-\lambda)\eta};q^{(2k-\lambda
)\eta})_{\infty}.
\end{align}

Multiplying  both sides of the above identity by
\begin{equation}\label{sim-4}
\frac{(-q^{\eta-\alpha_1},\ldots,-q^{\eta-\alpha_{\lambda}},
-q^\eta;q^\eta)_\infty}{(q^{\eta};q^{\eta})_{\infty}},
\end{equation}
and noting that for $1\leq i\leq\lambda$, $\alpha_i+\alpha_{\lambda+1-i}=\eta$, we show that Theorem  \ref{G-B-O-1-eq} holds when $k>r>\lambda\geq 0$.

Case 2. When $k=r>\lambda\geq 0$. We first determine the desired Bailey pair. To this end, we begin with  the unit Bailey pair \cite[H(17)]{Slater-1952},
\begin{align*}
 & \alpha_n^{(0)}(1, q)= \begin{cases}1, & \text { if } n=0, \\
  (-1)^n q^{n^2 / 2}\left(q^{-n / 2}+q^{n / 2}\right), & \text { if } n \geq 1,\end{cases} \\[5pt]
 & \beta_n^{(0)}(1, q)=\delta_{n,0}= \begin{cases}1, & \text { if } n=0, \\
  0, & \text { if } n \geq 1 .\end{cases}
\end{align*}

Applying Proposition 2.4 with $A={1}/{2}$ to ($\alpha_n^{(0)}(1, q),\,\beta_n^{(0)}(1, q)$), we get the following Bailey pair,
\begin{align} \label{baily-special}
& \alpha_n(1, q)= \begin{cases}1, & \text { if } n=0, \\
(-1)^n q^{n^2 / 2}\left(q^{-n / 2}+q^{-n / 2}\right)(1+q^n)/2, & \text { if } n \geq 1,\end{cases}\nonumber \\[5pt]
& \beta_n(1, q)= \frac{1+q^n}{2}\delta_{n,0}.
\end{align}

 Substitute this Bailey pair with $q$ replaced by $q^\eta$ into Proposition \ref{lc} with $r=k$.  Under the assumption  that  $\alpha_i+\alpha_{\lambda+1-i}=\eta$ for $1\leq i\leq \lambda$, and applying  Jacobi's triple product identity, we derive that the  left-hand side of \eqref{lem-cor} can be simplified as follows.
\begin{align}\nonumber
& 1+\sum_{n=1}^\infty\frac{(-q^{\alpha_1};q^{\eta})_n\cdots(-q^{
		\alpha_{\lambda}};q^{\eta})_n}{(-q^{\eta-\alpha_1};q^{\eta})
	_n\cdots(-q^{\eta-\alpha_{\lambda}};q^{\eta})_n}\\[5pt] \nonumber
& \quad\times (-1)^nq^{(k-\frac{\lambda}{2})\eta n^2+\frac{\lambda+1}{2}\eta n
	-(\alpha_1+\cdots+\alpha_{\lambda})n}(q^{-\frac{1}{2}\eta n}+q^{-\frac{1}{2}\eta n})\\ \nonumber
&\qquad=1+2\sum_{n=1}^{\infty}(-1)^nq^{(k-\frac{\lambda}{2})\eta n^2}\\ \label{eq-L-2}
&\quad\qquad=(q^{(k-\frac{\lambda}{2})\eta},q^{(k-
	\frac{\lambda}{2})\eta},q^{(2k-\lambda)\eta};q^{(2k-\lambda)\eta})_{\infty}.
\end{align}

When applying $\beta_n(1,q^\eta)$ from \eqref{baily-special} to \eqref{lem-cor} with $r=k$, we discover that the terms with $N_k>0$ in the summation of the right-hand side of \eqref{lem-cor} equals zero. Therefore, the right-hand side of \eqref{lem-cor} should be
\begin{align} \nonumber
&\frac{(q^{\eta};q^{\eta})_{\infty}}{(-q^{\eta-\alpha_1};q^{\eta})_
{\infty}}\sum_{N_1\geq N_2\geq\cdots\geq N_{k-1}\geq 0}
\frac{q^{
\eta(N_{\lambda+2}^2+\cdots+N_{k-1}^2)+\eta\left({N_1+1 \choose2}+\cdots+{N_{\lambda+1}+1 \choose2}\right)
-(\alpha_1N_1+\cdots+\alpha_{\lambda}N_{\lambda})
}}{(q^{\eta};q^{\eta})_{N_1-N_{2}}\cdots(q^{\eta};q^{\eta})_{
N_{k-1}}}\\[5pt]
&\qquad\times\frac{(-1;q^{\eta})_{N_{\lambda+1}
}(-q^{\alpha_1};q^{\eta})_{N_1}\cdots(-q^{\alpha_{\lambda
}};q^{\eta})_{N_{\lambda}}}
{(-q^{\eta};q^{\eta})_{N_{\lambda}}
(-q^{\eta-\alpha_2};q^
{\eta})_{N_1}\cdots(-q^{\eta-\alpha_{\lambda}};q^{\eta})_{
N_{\lambda-1}}} \nonumber \\ \nonumber
&=\frac{(q^{\eta};q^{\eta})_{\infty}}{(-q^{\eta-\alpha_1};q^{\eta})_
	{\infty}}\sum_{N_1\geq\cdots\geq N_{k-1}\ge 0}\frac{2q^{\eta(N_1^2+\cdots+N_{k-1}^2)}
	(-q^{\eta-\eta N_{\lambda+1}};q^{\eta})_{N_{\lambda+1}-1}}
{(q^{\eta};q^{\eta})_{N_1-N_2}\cdots(q^{\eta};q^{\eta})_
	{N_{k-1}}}\\[5pt] \label{eq-R2-2}
&\hskip1cm\times\frac{(-q^{\eta+\eta N_{\lambda}};q^{\eta})_{\infty}\prod_{s=1}^{\lambda}(-q^{\eta-\alpha_s
		-\eta N_{s}};q^{\eta})_{N_s}\prod_{s=2}^{\lambda}(-q^{\eta-\alpha_s+
		\eta N_{s-1}};q^{\eta})_{\infty}}{(-q^\eta;q^\eta)_\infty\prod_{s=2}^{\lambda}(-q^{
		\eta-\alpha_s};q^{\eta})_{\infty}},
\end{align}
where the last equality follows from \eqref{sim-3}, \eqref{sim-1} and \eqref{sim-2}.

Combining \eqref{eq-L-2} and \eqref{eq-R2-2}, we have
\begin{align}\nonumber
&\frac{(q^{\eta};q^{\eta})_{\infty}}{(-q^{\eta-\alpha_1};q^{\eta})_
	{\infty}}\sum_{N_1\geq\cdots\geq N_{k-1}\geq0}\frac{2q^{\eta(N_1^2+\cdots+N_{k-1}^2)}
	(-q^{\eta-\eta N_{\lambda+1}};q^{\eta})_{N_{\lambda+1}-1}}
{(q^{\eta};q^{\eta})_{N_1-N_2}\cdots(q^{\eta};q^{\eta})_{N_{k-1}}}\\[5pt] \nonumber
&\quad\times\frac{(-q^{\eta+\eta N_{\lambda}};q^{\eta})_{\infty}\prod_{s=1}^{\lambda}(-q^{\eta-\alpha_s
		-\eta N_{s}};q^{\eta})_{N_s}\prod_{s=2}^{\lambda}(-q^{\eta-\alpha_s+
		\eta N_{s-1}};q^{\eta})_{\infty}}{(-q^\eta;q^\eta)_\infty\prod_{s=2}^{\lambda}(-q^{
		\eta-\alpha_s};q^{\eta})_{\infty}}\\ \nonumber
&\qquad=(q^{(k-\frac{\lambda}{2})\eta},q^{(k-
	\frac{\lambda}{2})\eta},q^{(2k-\lambda)\eta};q^{(2k-\lambda
	)\eta})_{\infty}.
\end{align}
By multiplying  both sides of the above identity by  \eqref{sim-4},
and using the assumption that $\alpha_i+\alpha_{\lambda+1-i}=\eta$ for $1\leq i\leq\lambda$,  we deduce that Theorem \ref{G-B-O-1} also holds when $k=r>\lambda\geq 0$. This completes the proof. \qed

\section{Proof of Theorem  \ref{Over-bre-1a}}

In this section, we first assert that it suffices to show the following result in order to prove Theorem \ref{Over-bre-1a}.

\begin{thm}\label{lambda-r} For   $k> r\geq \lambda\geq 2$,
	\begin{equation}\label{lambda}
	\begin{split} &\sum_{n\geq0}\overline{B}_1(\alpha_{1},\ldots,\alpha_{\lambda};\eta,k,r;n)q^n\\[5pt]
	&\quad=
	(-q^{\alpha_{1}},-q^{\alpha_{\lambda}};q^{\eta})_\infty
	 \sum_{n\geq0}\overline{B}_1(\alpha_{2},\ldots,\alpha_{\lambda-1};\eta,k-1,r-1;n)q^n.
	\end{split}
	\end{equation}
\end{thm}
Before proceeding to the proof of Theorem \ref{lambda-r},  we first demonstrate the derivation of Theorem \ref{Over-bre-1a} from Theorem \ref{lambda-r}.

 \subsection{Proof of Theorem \ref{Over-bre-1a} with the aid of Theorem  \ref{lambda-r}}

\noindent{\it Proof of Theorem \ref{Over-bre-1a}.}  By induction on $\lambda$. When $\lambda=0$,    setting $q\rightarrow q^\eta$ in the following generating function due to  Chen, Sang and Shi \cite[(1.1)]{Chen-Sang-Shi-2013}:
\begin{equation*}\label{proof-0-A-0}
\sum_{n\geq0}\overline{B}_1(-;1,k,r;n)q^n=\frac{{(-q ;q )_\infty}(q^{r},
	q^{2k-r},
	q^{2k};
	q^{2k})_\infty}{(q ;q )_\infty},
\end{equation*}
where   $k> r\geq1$,  we find that Theorem \ref{Over-bre-1a} holds when $\lambda=0$.

When $\lambda=1$, observing that $\eta=2\alpha_1$, we see that $\eta$
must be even  and $\alpha_1=\eta/2$. It follows from \cite[Theorem 1.19]{he-ji-zhao} that for $k>r\geq 1$,
 \begin{equation}\label{odd-2}
 \begin{split}
 &\sum_{n\geq0}\overline{B}_1(1;2,k,r;n)q^n= \frac{(-q;q^2)_\infty(-q^2;q^2)_\infty(q^{2r-1},q^{4k-2r-1},
 q^{4k-2};q^{4k-2})_\infty}{(q^2;q^2)_\infty}.
 \end{split}
 \end{equation}
Letting  $q\rightarrow q^{\eta/2}$ in \eqref{odd-2}, we deduce that Theorem \ref{Over-bre-1a}  holds when $\lambda=1$.

When $\lambda\geq 2$, assume that Theorem \ref{Over-bre-1a} holds for   $\lambda-2$. Then,  for  $k-1> r-1> \lambda-2\geq0$, we have
\begin{equation}\label{proof-0-A-i}
 \begin{split}
 &\sum_{n\geq0}\overline{B}_1(\alpha_2,\ldots,\alpha_{\lambda-1};\eta,k-1,r-1;n)q^n\\
 &\quad=
\frac{{(-q^{\alpha_2},\ldots,-q^{\alpha_{\lambda-1}},-q^\eta;q^\eta)_\infty}
(q^{\eta(r-\frac{\lambda}{2})},
 q^{\eta(2k-r-\frac{\lambda}{2})},
 q^{\eta(2k-\lambda)};
 q^{\eta(2k-\lambda)})_\infty}{(q^\eta;q^\eta)_\infty}.
 \end{split}
 \end{equation}
By substituting \eqref{proof-0-A-i} into \eqref{lambda}, we obtain \eqref{proof-0-A}. This implies that Theorem \ref{Over-bre-1a} holds for $\lambda$. Thus, we  have proven Theorem \ref{Over-bre-1a} with the aid of Theorem \ref{lambda-r}. \qed

\subsection{The outline of the proof of Theorem \ref{lambda-r} }

Let $\mathcal{D}_{\alpha_1}$ and $\mathcal{D}_{\alpha_\lambda}$ denote the sets of distinct partitions whose parts are congruent to $\alpha_1$ and $\alpha_\lambda $ modulo $\eta$,  respectively.  Clearly, we have
\[\sum_{\delta\in\mathcal{D}_{\alpha_1}}q^{|\delta|}
=(-q^{\alpha_1};q^{\eta})_\infty\quad\text{and}\quad\sum_{\delta\in
\mathcal{D}_{\alpha_\lambda}}q^{|\delta|}=(-q^{\alpha_\lambda};q^{\eta})_\infty.\]
 For $k> r\geq\lambda\geq0$, let $\overline{\mathcal{B}}_1(\alpha_1,\ldots,\alpha_\lambda;\eta,k,r)$ denote the set of partitions  counted by $\overline{B}_1(\alpha_1,\ldots,\alpha_\lambda;\eta,k,r;n)$ for $n\geq 0$. It is easy to see that  Theorem \ref{lambda-r} is equivalent to the following combinatorial statement.

\begin{thm}\label{lambdathm} For $k>r\geq\lambda\geq2$, there is a bijection $\Theta$ between  $\mathcal{D}_{\alpha_1}\times\mathcal{D}_{\alpha_\lambda}\times \mathcal{\overline{B}}_1(\alpha_{2},\ldots,\alpha_{\lambda-1};\eta,
k-1,r-1)$ and $\mathcal{\overline{B}}_1(\alpha_1, \ldots,\alpha_{\lambda};\eta,k,r)$. Moreover, for   a triple $(\delta^{(1)},\delta^{(\lambda)},\pi) \break $$ \in
\mathcal{D}_{\alpha_1}\times\mathcal{D}_{\alpha_\lambda}\times \mathcal{\overline{B}}_1(\alpha_{2},\ldots,\alpha_{\lambda-1} ;\eta,
k-1,r-1)$, we have $\tau=\Theta(\delta^{(1)},\delta^{(\lambda)},\pi)\in \mathcal{\overline{B}}_1(\alpha_1,\ldots,\alpha_{\lambda};\eta, \break $$ k,r)$ such that  $|\tau|=|\delta^{(1)}|+|\delta^{(\lambda)}|+|\pi|$.
\end{thm}

To prove this theorem, we adopt the analogous strategy used for Conjecture \ref{Bressoud-conjecture-j} for $j=1$ given by Kim \cite{Kim-2018}. Let $\mathcal{\overline{{C}}}_1(\alpha_2,\ldots,\alpha_\lambda;\eta,k,r)$ denote the set of overpartitions in $\mathcal{\overline{{B}}}_1(\alpha_1,\ldots,\alpha_\lambda;\eta,k,r)$ without parts $\equiv \alpha_1\pmod{\eta}$. Bear in mind that $\alpha_2,\ldots,\alpha_\lambda$ and $\eta$ are integers such that
\[
0<\alpha_2<\cdots<\alpha_\lambda<\eta, \quad \text{and} \quad \alpha_i=\eta-\alpha_{\lambda+1-i}\quad \text{for} \quad 2\leq i\leq \lambda-1.
\]
To build the bijection $\Theta$ in Theorem \ref{lambdathm}, we will first unite $\pi$ and $\delta^{(\lambda)}$, and denote the resulting overpartition by $\pi^{(0)}$. Evidently, $\pi^{(0)} \in \mathcal{\overline{{C}}}_1(\alpha_2,\ldots,\alpha_\lambda;\eta,k,r)$.
We next  aim to merge the parts of  $\delta^{(1)}$ into $\pi^{(0)}$ from smallest to largest to generate an overpartition in $\mathcal{\overline{B}}_1(\alpha_1,\ldots,\alpha_{\lambda};\eta,k,r)$.  There are two steps.   We first   merge some parts $\equiv \alpha_1\pmod{\eta}$ in $\delta^{(1)}$ and some parts $\equiv \alpha_\lambda \pmod{\eta}$ in $\pi^{(0)}$ to generate some non-overlined parts divisible by $\eta$ (due to the fact that $\alpha_1+\alpha_\lambda=\eta$). It turns out the $(k-1)$-addition  introduced by Kim  \cite{Kim-2018} can fulfill this task directly (see Section 3.4). In the second step,  we are   meant to merge
the remaining parts of $\delta^{(1)}$ (which are $\equiv \alpha_1\pmod{\eta}$) and the   overpartition in $\mathcal{\overline{{C}}}_1(\alpha_2,\ldots,\alpha_\lambda;\eta,k,r)$ to generate certain overlined parts $\equiv \alpha_1\pmod{\eta}$. As a result, we get an overpartition in $\mathcal{\overline{{B}}}_1(\alpha_1,\ldots,\alpha_\lambda;\eta,k,r)$. We can accomplish our objective by utilizing the $(k-1)$-insertion outlined in \cite{he-ji-zhao} with $a=\alpha_1$ (see Section 3.5).

In the next subsection,  we will first recall the definitions of the Gordon marking and the  reverse Gordon marking of an overpartition in $\overline{\mathcal{B}}_1(\alpha_1,\ldots,\alpha_\lambda;\eta,k,r)$ established in our previous paper \cite{he-ji-zhao}. The $(k-1)$-addition operation and the $(k-1)$-insertion operation are defined based on the reverse Gordon marking, while their inverse operations are defined based on the  Gordon marking. We also review the forward move and the backward move, which are the main ingredients in the constructions of the $(k-1)$-addition, the $(k-1)$-insertion and their inverse operations. In Section 3.4, we recall the $(k-1)$-addition and its inverse map (i.e., the $(k-1)$-subtraction) given by Kim \cite{Kim-2018} in the context of overpartitions. Section 3.5 is devoted to revisiting the $(k-1)$-insertion and its inverse map (i.e., the $(k-1)$-separation) defined in \cite{he-ji-zhao}. In Section 3.6, we give a proof of Theorem \ref{lambdathm}  by successively applying the $(k-1)$-addition operation and the $(k-1)$-insertion operation.

\subsection{The (reverse) Gordon marking and the forward (backward) move}

This subsection  revisits  the definitions of the Gordon marking and the  reverse Gordon marking of an overpartition in $\overline{\mathcal{B}}_1(\alpha_1,\ldots,
\alpha_\lambda;\eta,k,r)$,  originally established in our earlier paper \cite{he-ji-zhao}. We follow the terminology, notation and examples in \cite{he-ji-zhao}.

In the rest of this paper, we assume that $k$, $r$ and $\lambda$ are integers such that  $k> r\geq\lambda\geq2$. The Gordon marking of an overpartition in $\overline{\mathcal{B}}_1(\alpha_1,\ldots,\alpha_\lambda;\eta,k,r)$ was defined as follows.

\begin{defi}{\rm\cite[Definition 3.1]{he-ji-zhao}}\label{Gordon-marking}
Let $\pi=(\pi_1,\pi_2,\ldots,\pi_\ell)$ be an overpartition satisfying {\rm{(1)}} and {\rm{(2)}} in Definition \ref{defi-O-B}. Assign a positive integer to each part of $\pi$ as follows{\rm{:}} First,  assign $1$ to  $\pi_\ell$. Then, for each $\pi_i$, assign $s$  to $\pi_i$, where $s$ is the smallest positive integer that is not used to mark the parts  $\pi_m$ such that $m>i$ and $\pi_m\geq \pi_i-\eta$ with strict inequality if $\pi_i$ is overlined. Denote the Gordon marking of $\pi$ by $G(\pi)$.
\end{defi}

For example, let  $\pi$ be an overpartition in $\overline{\mathcal{B}}_1(1,5,9;10,5,4)$  given by
\begin{equation}\label{mark-exa-1}
\begin{split}
&\pi=(\overline{80},80,80,\overline{70},70,\overline{69},\overline{60},
60,\overline{55},\overline{51},{50},\overline{49},\overline{45},\overline{41},\overline{39},\overline{35},\\[3pt]
&\ \ \ \ \ \ \ \overline{29},\overline{20},{20},20,\overline{11},\overline{10},\overline{9},
\overline{5},\overline{1}).
\end{split}
\end{equation}
The Gordon marking of $\pi$ is given by
\begin{equation*}\label{exa-g-1}
\begin{split}
&G(\pi)=(\overline{80}_1,{ {80}_4},{80}_2,\overline{70}_1,{70}_3,\overline{69}_2,{ \overline{60}_4},
{60}_1,\overline{55}_2,\overline{51}_3,{ {50}_4},\overline{49}_1,\overline{45}_2,\overline{41}_3,\overline{39}_1,\overline{35}_2,\\[3pt]
&\ \ \ \ \ \ \ \ \ \ \ \  \overline{29}_1,{ \overline{20}_4},{20}_3,{20}_2,\overline{11}_1,{ \overline{10}_4},\overline{9}_3,
\overline{5}_2,\overline{1}_1),
\end{split}
\end{equation*}
where the subscript of each part  represents the mark in the Gordon marking.

Letting $\pi$ be an overpartition satisfying {\rm(1)} and {\rm(2)} in Definition \ref{defi-O-B}, by definition, we see that  $\pi$  is an overpartition in $\overline{\mathcal{B}}_1(\alpha_1,\ldots,\alpha_\lambda;\eta,k,r)$ if and only if the marks in $G(\pi)$ do not exceed  $k-1$ and  the marks of parts less than or equal to $\eta$ in $G(\pi)$ do not exceed    $r-1$.

The reverse Gordon marking of an overpartition in $\overline{\mathcal{B}}_1(\alpha_1,\ldots,\alpha_\lambda;\eta,k,r)$  is defined by assigning a mark to each part starting with the largest part instead.

\begin{defi}{\rm\cite[Definition 3.2]{he-ji-zhao}}\label{R-Gordon-marking}
Let $\pi=(\pi_1,\pi_2,\ldots,\pi_\ell)$ be an overpartition satisfying {\rm{(1)}} and {\rm{(2)}} in Definition \ref{defi-O-B}. Assign a positive integer to each part of $\pi$ as follows{\rm{:}} First assign $1$ to $\pi_1$. Then, for each $\pi_i$, assign $s$ to $\pi_i$, where $s$ is   the smallest positive integer that is not used to mark the parts  $\pi_m$  such that $m<i$ and $\pi_m\leq \pi_i+\eta$ with strict inequality if $\pi_i$ is overlined. Denote the reverse Gordon marking of $\pi$ by $RG(\pi)$.
\end{defi}

For the  overpartition $\pi$ in $\overline{\mathcal{B}}_1(1,5,9;10,5,4)$  defined in \eqref{mark-exa-1}, the reverse Gordon marking  of $\pi$ reads
\begin{equation}\label{exa-r-1}
\begin{split}
&RG(\pi)=(\overline{80}_1,{80}_2,{80}_3,\overline{70}_1,{ {70}_4},\overline{69}_2,\overline{60}_1,
{60}_3,\overline{55}_2,{ \overline{51}_4},{50}_1,\overline{49}_3,\overline{45}_2,{ \overline{41}_4},\overline{39}_1,\overline{35}_2,\\[3pt]
&\ \ \ \ \ \ \ \ \ \ \ \ \ \  \overline{29}_1,\overline{20}_2,{20}_3,{ {20}_4},\overline{11}_1,\overline{10}_2,\overline{9}_3,{ \overline{5}_4},\overline{1}_1).
\end{split}
\end{equation}

Analogously,   an overpartition $\pi$  satisfying {\rm(1)} and {\rm(2)} in Definition \ref{defi-O-B} is an overpartition in $\overline{\mathcal{B}}_1(\alpha_1,\ldots,\alpha_\lambda;\eta,k,r)$  if and only if the marks in $RG(\pi)$ do not exceed  $k-1$ and there are at most $r-1$ parts less than or equal to $\eta$ in $\pi$.

We proceed to recall the definition of $(k-1)$-bands. Let  $\pi=(\pi_1,\pi_2,\ldots,\pi_\ell)$ be
 an overpartition in $\overline{\mathcal{B}}_1(\alpha_1,\ldots,\alpha_\lambda;\eta,k,r)$. If there are  $k-1$ consecutive parts $\pi_i\geq \pi_{i+1}\geq \cdots \geq \pi_{i+k-2}$ satisfying the following relation:
 \begin{equation*}\label{sequence}
\pi_i\leq\pi_{i+k-2}+\eta\text{ with strict inequality if }
\\ \pi_i \text{ is overlined},
\end{equation*}
 then such $k-1$ parts will be called a $(k-1)$-band of $\pi$.
Observe that for a $(k-1)$-band $\{\pi_{i+l}\}_{0\leq l\leq k-2}$ of $\pi$ without  overlined parts divisible by $\eta$, if we change
the overlined parts of $\pi$ to non-overlined parts, then this $(k-1)$-band reduces to a $(k-1)$-sequence introduced by Kim \cite{Kim-2018}.

For example,  let $\pi$ be the overpartition in $\overline{\mathcal{B}}_1(1,5,9;10,5,4)$ defined in \eqref{mark-exa-1}, where $k=5$. There are twelve $4$-bands in  $\pi$.
  \[\{80,{ 80},\overline{70},{ 70}\},\{{ 70},\overline{69},{ \overline{60}},60\},
  \{{ \overline{60}},{60},\overline{55},{ \overline{51}}\},
  \{{60},\overline{55},{ \overline{51}},{ {50}}\},\]
   \[\{\overline{55},{ \overline{51}},{ {50}},\overline{49}\},
   \{{ \overline{51}},{ {50}},\overline{49},\overline{45}\},\{{{50}},\overline{49},
   \overline{45},{ \overline{41}}\},\{\overline{29},{ \overline{20}},{20},{ 20}\},\]
   \[\{{ \overline{20}},{20},{ 20},\overline{11}\},\{20,{ 20},\overline{11},{ \overline{10}}\},
   \{\overline{11},{ \overline{10}},\overline{9},{ \overline{5}}\},\{{ \overline{10}},\overline{9},{ \overline{5}},\overline{1}\}.\]

For each $(k-1)$-band $\{\pi_{i+l}\}_{0\leq l\leq k-2}$ of $\pi$,
it is easy to see that the marks of $\pi_{i+l}$ are distinct in
the Gordon marking and the reverse Gordon marking of  $\pi$.
Hence there exists one part in  $\{\pi_{i+l}\}_{0\leq l\leq k-2}$ marked with $k-1$ in the Gordon marking and the reverse Gordon marking of  $\pi$.  We  now restrict our attention to two kinds of special $(k-1)$-bands.

{\it The $(k-1)$-bands of the first kind} refer to those bands in which the $(k-1)$-marked part in the Gordon marking is the largest element.
 Assume that there are $N$ parts marked with $k-1$ in $G(\pi)$,
and denote these $(k-1)$-marked parts by  $\tilde{g}_1(\pi)> {\tilde{g}}_2(\pi)>\cdots >\tilde{g}_N(\pi)$. For each $(k-1)$-marked part $\tilde{g}_p(\pi)$
in $G(\pi)$,    there is a $(k-1)$-band of $\pi$ such that $\tilde{g}_p(\pi)$ is the largest element of this $(k-1)$-band.
Such a $(k-1)$-band is called the $(k-1)$-band
induced by $\tilde{g}_p(\pi)$, denoted  $\{\tilde{g}_p(\pi)\}_{k-1}$.

For example, for   the  overpartition   $\pi$ given in \eqref{mark-exa-1},
  there are five $4$-marked parts in $G(\pi)$, namely,
 $\tilde{g}_1(\pi)=80$, $\tilde{g}_2(\pi)=\overline{60}$,
 $\tilde{g}_3(\pi)={50}$, $\tilde{g}_4(\pi)=\overline{20}$ and $\tilde{g}_5(\pi)=\overline{10}$.  The   $4$-bands induced by  $\tilde{g}_1(\pi),\tilde{g}_2(\pi),\tilde{g}_3(\pi),
 \tilde{g}_4(\pi)$  and $\tilde{g}_5(\pi)$   are  illustrated
 in $G(\pi)$ below:
\begin{equation*} \label{mark-exa-1-s}
\begin{split}
&G(\pi)=(\overline{80}_1,\overbrace{{ {80}_4},
{ {80}_2,
\overline{70}_1,{70}_3}}^{{\{80\}_4}},\overline{69}_2,\overbrace{{ \overline{60}_4},{ {60}_1,\overline{55}_2,\overline{51}_3}}^{{\{\overline{60}\}_4}},\overbrace{{ {50}_4},{ \overline{49}_1,\overline{45}_2,\overline{41}_3}}^{{\{{50}\}_4}},\overline{39}_1,\overline{35}_2,\\[5pt]
&\ \ \ \ \ \ \ \ \ \ \
\overline{29}_1,\underbrace{{ \overline{20}_4},{ {20}_3,
{20}_2,\overline{11}_1}}_
{\{\overline{20}\}_4},\underbrace{{ {\overline{10}_4}},
{ \overline{9}_3,\overline{5}_2,\overline{1}_1}}_{\{\overline{10}\}_4}).
\end{split}
\end{equation*}

{\it The $(k-1)$-bands of the second kind} are a specific category of bands that are characterized by the $(k-1)$-marked part in the reverse Gordon marking being
the smallest part. Assume that there are $M$ parts marked with $k-1$  in $RG(\pi)$, namely,  $\tilde{r}_1(\pi)> \tilde{{r}}_2(\pi)>\cdots >\tilde{r}_M(\pi)$.
By the same reasoning, we see that  there is a $(k-1)$-band of $\pi$ in which $\tilde{r}_p(\pi)$ is the smallest element.  Such a $(k-1)$-band is called the $(k-1)$-band induced by $\tilde{r}_p(\pi)$, denoted  $\{\tilde{r}_p(\pi)\}_{k-1}$.

For example, for the overpartition  $\pi$    given in \eqref{mark-exa-1},   there are five $4$-marked parts in $RG(\pi)$, which are  $\tilde{r}_1(\pi)=70$, $\tilde{r}_2(\pi)=\overline{51}$, $\tilde{r}_3(\pi)=\overline{41}$, $\tilde{r}_4(\pi)={20}$  and $\tilde{r}_5(\pi)=\overline{5}$. The   $4$-bands induced by  $\tilde{r}_1(\pi)$, $\tilde{r}_2(\pi)$, $\tilde{r}_3(\pi)$, $\tilde{r}_4(\pi)$ and  $\tilde{r}_5(\pi)$ are  displayed  below:
\begin{equation*} \label{rmark-exa-1-s}
\begin{split}
&RG(\pi)=(\overline{80}_1,\overbrace{{ {80}_2,{80}_3,\overline{70}_1},
{ {70}_4}}^{{\{70\}_4}},\overline{69}_2,
\overbrace{{ \overline{60}_1,{60}_3,\overline{55}_2}, { \overline{51}_4}}^{{\{\overline{51}\}_4}},\overbrace{{{50}_1,\overline{49}_3,\overline{45}_2}, { \overline{41}_4}}^{{\{\overline{41}\}_4}}, \overline{39}_1,\overline{35}_2,\\[5pt]
&\ \ \ \ \ \ \ \ \ \ \ \ \ \ \underbrace{{ \overline{29}_1,\overline{20}_2,} { {20}_3},{ {20}_4}}
_{\{20\}_4},
\underbrace{{ \overline{11}_1,\overline{10}_2,\overline{9}_3},
{ \overline{5}_4}}_{\{\overline{5}\}_4},\overline{1}_1).
\end{split}
\end{equation*}

The following proposition indicates that the number of $(k-1)$-marked parts in $G(\pi)$ equals the number of $(k-1)$-marked parts in  $RG(\pi)$.
\begin{prop}{\rm\!\! \cite[Proposition 3.3]{he-ji-zhao}} \label{sequence-length}
Let $\pi$ be an overpartition in $\overline{\mathcal{B}}_1
 (\alpha_1,\ldots,\alpha_\lambda;\eta,k,r)$. Assume that there are $N$ parts marked with $k-1$ in $G(\pi)$, say, $\tilde{g}_1(\pi)> \tilde{g}_2(\pi)>\cdots >\tilde{g}_N(\pi)$, and there are  $M$ parts marked with $k-1$ in $RG(\pi)$, say, $\tilde{r}_1(\pi)> \tilde{r}_2(\pi)>\cdots >\tilde{r}_M(\pi)$. Then $N=M$.   Moveover, for each $1\leq i \leq N$, we have $\tilde{g}_i(\pi)\in \{\tilde{r}_i(\pi)\}_{k-1}$ and $\tilde{r}_i (\pi)\in \{\tilde{g}_i(\pi)\}_{k-1}$, where $\{\tilde{g}_i(\pi)\}_{k-1}$ {\rm(}resp. $\{\tilde{r}_i(\pi)\}_{k-1}${\rm)} is the $(k-1)$-band of $\pi$ induced by $\tilde{g}_i(\pi)$ {\rm(}resp. $\tilde{r}_i(\pi)${\rm)}.
\end{prop}

We conclude this subsection by stating  the forward move and  the backward move  based on the Gordon marking and the reverse Gordon marking of an overpartition in $\overline{\mathcal{B}}_1(\alpha_1,\ldots,\alpha_\lambda;\eta,k,r)$, respectively. For more details about them, please refer to Section 4.1 in  \cite{he-ji-zhao}.

\begin{defi}{\rm\cite[Definition 4.2]{he-ji-zhao}}\label{forward-defi} For  $N\geq 1$, let $\pi$ be an overpartition satisfying {\rm(1)}, {\rm(2)} and  {\rm(3)}  in Definition \ref{defi-O-B}. Assume that there are $N$ parts marked with $k-1$ in  $RG(\pi)$, say  $\tilde{r}_1(\pi)>\tilde{r}_2(\pi)>\cdots>
\tilde{r}_N(\pi)$.  For $1\leq p\leq N$,  the forward  move  $\phi_p$ is defined as follows{\rm :} add $\eta$ to each of  $\tilde{r}_1(\pi),\,\tilde{r}_2(\pi),\ldots,\tilde{r}_{p}(\pi)$ and  rearrange the parts in  non-increasing order to obtain a new overpartition, denoted $\phi_p(\pi)$.
\end{defi}

When $p=0$, the forward move $\phi_p$ is defined to be the identity map, that is, $\phi_p(\pi)=\pi$.

For example, let $\pi$ be  the overpartition defined in  \eqref{mark-exa-1}, whose reverse Gordon marking is given in \eqref{exa-r-1}. Apply the forward move $\phi_3$ to  $\pi$, namely, add $\eta=10$ to each of $\tilde{r}_1(\pi)=70$, $\tilde{r}_2(\pi)=\overline{51}$ and $\tilde{r}_3(\pi)=\overline{41}$, and  so we get
 \begin{equation*}\label{exa-r-2}
 \begin{split}
&\phi_3(\pi)=(\overline{80},{ {80}},{80},{80},\overline{70},\overline{69},{\overline{61}},\overline{60},
{60},\overline{55},{\overline{51}},{50},\overline{49},\overline{45},\overline{39},\overline{35},\\
&\ \ \ \ \ \ \ \ \ \ \ \ \ \overline{29},\overline{20},{20},{20},\overline{11},
\overline{10},\overline{9},\overline{5},\overline{1}).
\end{split}
\end{equation*}

\begin{defi}{\rm\cite[Definition 4.4]{he-ji-zhao}}\label{defi-backward}
For $N\geq p\geq 1$, let $\omega$ be an overpartition satisfying {\rm(1)}, {\rm(2)} and  {\rm(3)}  in Definition \ref{defi-O-B}.  Assume that there are $N$ parts marked with $k-1$ in  $G(\omega)$, denoted   $\tilde{g}_1(\omega)>\tilde{g}_2(\omega)>\cdots>
\tilde{g}_N(\omega)$, for which $\tilde{g}_p(\omega)\geq\overline{\eta+\alpha_1}$. The backward  move $\psi_p$ is defined as follows{\rm :} subtract $\eta$ from each  of $\tilde{g}_1(\omega),\,\tilde{g}_2(\omega),\ldots,
 \tilde{g}_{p}(\omega)$ and  rearrange the parts in  non-increasing order to obtain a new overpartition, denoted  $\psi_p(\omega)$.
\end{defi}

For example, let $\pi$ be the overpartition in $\overline{\mathcal{B}}_1(1,5,9;10,5,4)$ given in    \eqref{mark-exa-1},  and let $\omega=\phi_3(\pi)$. Then
the Gordon marking of $\omega$ is given by
\begin{equation*}\label{example-b}
\begin{split}
&G(\omega)=(\overline{80}_1,\overbrace{{ {80}_4},
{ {80}_3,
{80}_2,\overline{70}_1}}^{\{80\}_4},\overline{69}_2,\overbrace{{ \overline{61}_4},{ \overline{60}_3,{60}_1,\overline{55}_2}}^{\{\overline{61}\}_4},\overbrace{{ \overline{51}_4},{ {50}_3,\overline{49}_1,\overline{45}_2}}^{\{\overline{51}\}_4},\overline{39}_1,\overline{35}_2,\\[5pt]
&\ \ \ \ \ \ \ \ \ \ \ \
\overline{29}_1,\underbrace{{ \overline{20}_4},{ {20}_3,
{20}_2,\overline{11}_1}}_{\{\overline{20}\}_4},\underbrace{{ \overline{10}_4},
{ \overline{9}_3,\overline{5}_2,\overline{1}_1}}_{\{\overline{10}\}_4}).
\end{split}
\end{equation*}
 There are five $4$-marked parts in $G(\omega)$, which are $\tilde{g}_1(\omega)=80$, $\tilde{g}_2(\omega)=\overline{61}$,
 $\tilde{g}_3(\omega)=\overline{51}$, $\tilde{g}_4(\omega)=\overline{20}$ and $\tilde{g}_5(\omega)=\overline{10}$. The backward move  $\psi_3$ transforms  $\omega$ back to $\pi$. That is, the overpartition $\pi$ can be obtained from $\omega$ by subtracting $\eta=10$ from each of $\tilde{g}_1(\omega)=80$, $\tilde{g}_2(\omega)=\overline{61}$ and
 $\tilde{g}_3(\omega)=\overline{51}$.

\subsection{The $(k-1)$-addition and the $(k-1)$-subtraction}

 Just as  mentioned before, to build the bijection $\Theta$ in Theorem \ref{lambdathm}, we first aim to  merge some parts $\equiv \alpha_1\pmod{\eta}$ in $\delta^{(1)}$ and some parts $\equiv \alpha_\lambda \pmod{\eta}$ in $\pi^{(0)}$ to generate some non-overlined parts divisible by $\eta$ (due to the fact that $\alpha_1+\alpha_\lambda=\eta$).  We find that the $(k-1)$-addition operation and the $(k-1)$-subtraction operation introduced by Kim \cite{Kim-2018} are capable of achieving this objective. Here we will adapt the $(k-1)$-addition operation and the $(k-1)$-subtraction operation introduced by Kim \cite{Kim-2018} for ordinary partitions to the setting of overpartitions in $\mathcal{\overline{{C}}}_1(\alpha_2,\ldots,\alpha_\lambda;\eta,k,r)$.

To present the definitions of the $(k-1)$-addition and the $(k-1)$-subtraction, we will define the non-degenerate $(k-1)$-bands of  an overpartition  $\pi \in \mathcal{\overline{{C}}}_1(\alpha_2,\ldots,
\alpha_\lambda;\eta,k,r)$ and the non-degenerate parts of an overpartition $\pi \in \mathcal{\overline{{C}}}_1(\alpha_2,\ldots,
\alpha_\lambda;\eta,k,r)$.

A $(k-1)$-band of $\pi$ is called a {\it  non-degenerate $(k-1)$-band}  if there are no parts $\equiv \alpha_\lambda \pmod{\eta}$ in this $(k-1)$-band.

To define the non-degenerate parts of an overpartition  $\pi \in \mathcal{\overline{{C}}}_1(\alpha_2,\ldots,
\alpha_\lambda;\eta,k,r)$, we first need to define the non-degenerate $(k-1)$-parts and the non-degenerate $(r-1)$-part.
For an overpartition $\pi$, we use  $f_\pi(0,\eta]$ to denote the number of parts  of $\pi$ less than or equal to $\eta$.

Let $\{\pi_{m+l}\}_{0\leq l\leq k-2}$ be a non-degenerate $(k-1)$-band  of $\pi$, namely,
\[\pi_m\geq \pi_{m+1}\geq \cdots \geq \pi_{m+k-2},\]
where $\pi_m\leq \pi_{m+k-2}+\eta$   with strict inequality  if  $\pi_m$   is overlined.
Note that $k>\lambda$, so  there is at least one non-overlined part in  $\{\pi_{m+l}\}_{0\leq l\leq k-2}$. Let $\pi_{m+t}$ be the largest non-overlined part in $\{\pi_{m+l}\}_{0\leq l\leq k-2}$. If $\pi_{m+t}>\pi_{m+t+1}$, then we call $\pi_{m+t}$ a {\it non-degenerate $(k-1)$-part of $\pi$}.

If $f_\pi(0,\eta]=r-1$ and $\overline{\alpha_\lambda}$ does not occur in $\pi$, then there is at least one non-overlined part $\eta$  of $\pi$  since $r-1>\lambda-2$.  Assume that  $\pi_{t}=\eta$, if $\pi_{t}>\pi_{t+1}$ and there are no $(k-1)$-bands  of $\pi$ in $(0,\overline{\eta+\alpha_\lambda})$,   then $\pi_{t}$ is called the {\it non-degenerate $(r-1)$-part of $\pi$}.   A part of $\pi$ is called a {\it non-degenerate  part} if it is either  a non-degenerate $(r-1)$-part of $\pi$ or a non-degenerate  $(k-1)$-part of $\pi$.

For example, let
  \begin{equation}\label{new-example-new-1}\begin{array}{lllllllllllllllll}
  &\ \pi_1,& \pi_2,& \pi_3,& \pi_4,& \pi_5,& \pi_6,& \pi_7,& \pi_8,& \pi_9,& \pi_{10},& \pi_{11}\\[5pt]
  &\ \downarrow&\downarrow&\downarrow&\downarrow&\downarrow&\downarrow
  &\downarrow&\downarrow&\downarrow&\downarrow&\downarrow\\[5pt]
 \pi=&(\overline{50},&\overline{47},&40,&40,&\overline{30},&\overline{25},&\overline{20},
&20,&\overline{15},&10,&\overline{5})
\end{array}\end{equation}
be the overpartition in $\mathcal{\overline{{C}}}_1(5,7;10,4,3)$, where $\eta=10,\,k=4,\,r=3,\,\lambda=3,\,\alpha_1=3,\, \alpha_2=5$ and $\alpha_3=7$. Note that $f_\pi(0,10]=2$,  there are no $3$-bands  of $\pi$ in $(0,\overline{17})$, and $\overline{7}$ does not occur in $\pi$. We see that $\pi_{10}=10$   is  the non-degenerate $2$-part of $\pi$.

For another example, we consider the following overpartition $\pi$ in $\mathcal{\overline{{C}}}_1(5,7;10,4,3)$, where $\eta=10,\,k=4,\,r=3,\,\lambda=3,\,\alpha_1=3,\, \alpha_2=5$ and $\alpha_3=7$.
  \begin{equation}\label{new-example-11}\begin{array}{lllllllllllllllll}
  &\ \pi_1,& \pi_2,& \pi_3,& \pi_4,& \pi_5,& \pi_6,& \pi_7,& \pi_8,& \pi_9,& \pi_{10},& \pi_{11}\\[5pt]
  &\ \downarrow&\downarrow&\downarrow&\downarrow&\downarrow&\downarrow
  &\downarrow&\downarrow&\downarrow&\downarrow&\downarrow\\[5pt]
 \pi=&(\overline{50},&\overline{47},&40,&40,&\overline{30},&\overline{25},&\overline{20},
&20,&\overline{10},&10,&\overline{7}).
\end{array}\end{equation}
It is easy to check that  there are three non-degenerate $3$-bands of $\pi$, which are $\{40,40,\overline{30}\}$, $\{\overline{25},\overline{20},20\}$ and $\{20,\overline{10},{10}\}$. Clearly, $\pi_4=40$ and $\pi_8=20$ are two non-degenerate $3$-parts of $\pi$.

 Here and in the sequel,
 we  make the following assumption. Let $\pi_i$ be the $i$-th part of the overpartition $\pi=(\pi_1,\pi_2,\ldots, \pi_\ell)$. If $\pi_i$ is an overlined part $\equiv \alpha_\lambda \pmod{\eta}$, then we define a new part $\pi_i+\alpha_1$ as a non-overlined part   of size $|\pi_i|+\alpha_1$. If $\pi_i$ is a non-overlined part divisible by $\eta$, then we define a new part $\pi_i-\alpha_1$ as an  overlined part   of size $|\pi_i|-\alpha_1$.

We will be concerned with the following two subsets of $\mathcal{\overline{{C}}}_1(\alpha_2,\ldots,\alpha_\lambda;\eta,k,r)$.

\begin{itemize}
\item For $0\leq p\leq N$,   let $\overline{\mathcal{C}}_\lambda(\alpha_2,\ldots,\alpha_\lambda;
\eta,k,r|N,p)$ denote the set of overpartitions $\pi$ in $\mathcal{\overline{{C}}}_1(\alpha_2,\ldots,\alpha_\lambda;\eta,k,r)$  such that   there are $N$ parts marked with $k-1$ in $RG(\pi)$, denoted $\tilde{r}_{1}(\pi)>\cdots >\tilde{r}_{N}(\pi)$, satisfying one of the following conditions:

\begin{itemize}
\item[(1)] If $0\leq p <N$, then there exists a part $\equiv \alpha_\lambda \pmod{\eta}$ in the $(k-1)$-band $\{\tilde{r}_{p+1}(\pi)\}_{k-1}$,  denoted $\tilde{\tilde{r}}_{p+1}(\pi)$, and there is no non-degenerate part of $\pi$ less  than $\tilde{\tilde{r}}_{p+1}(\pi)$;

\item[(2)] If $p=N$, then $f_\pi(0,\eta]=r-1$, $\overline{\alpha_\lambda}$ is a part of $\pi$, and $\tilde{r}_N(\pi)>\eta$ when $N\geq 1$.

\end{itemize}

\item For $0\leq p\leq N$,  let $\overline{\mathcal{C}}_\eta(\alpha_2,\ldots,\alpha_\lambda;\eta,k,r|N,p)$ denote the set of    overpartitions $\pi$ in $\mathcal{\overline{{C}}}_1(\alpha_2,\ldots,\alpha_\lambda;\eta,k,r)$  subject to the following conditions:

\begin{itemize}
\item[(1)] There exists a non-degenerate part of $\pi$;

\item[(2)] Let $\overline{\pi}$ be the overpartition obtained  from $\pi$ by subtracting $\alpha_1$ from the smallest non-degenerate part of $\pi$. Then there are $N$ parts marked with $k-1$ in the  Gordon marking  of $\overline{\pi}$, denoted  $\tilde{g}_{1}(\overline{\pi})>\cdots >\tilde{g}_{N}(\overline{\pi})$;

\item[(3)] If $0\leq p<N$, then there is no non-degenerate $(r-1)$-part of $\pi$ and   $p$ is the largest integer such that  $\tilde{g}_{p}(\overline{\pi})> \pi_m+\eta$, where $m$ is the largest integer such that  $\{\pi_{m+l}\}_{0\leq l\leq k-2}$ is a non-degenerate $(k-1)$-band of $\pi$;

\item[(4)] If $p=N$, then there is a non-degenerate $(r-1)$-part of $\pi$.

\end{itemize}

\end{itemize}

For example, let  $\pi$ be an overpartition in  $\mathcal{\overline{{C}}}_1(5,7;10,4,3)$, where $\eta=10,\,k=4,\,r=3,\,\lambda=3,\,\alpha_1=3,\, \alpha_2=5$ and $\alpha_3=7$. The reverse Gordon marking of $\pi$ is given below.
\begin{equation}\label{new-example-21}
RG(\pi)=(\overline{50}_1,\overline{47}_2,\overbrace{{{40}_1,\overline{30}_2},{{30}_3}}^{{\{{30}\}_3}},\overbrace{{ \overline{25}_1,\overline{20}_2},{ \overline{17}_3}}^{{ \{\overline{17}\}_3}},\overbrace{{ \overline{10}_1,{10}_2},{ \overline{7}_3}}^{{ \{\overline{7}\}_3}}).
\end{equation}
There are three parts marked with $3$ in $RG(\pi)$, which are $\tilde{r}_1(\pi)=30$, $\tilde{r}_2(\pi)=\overline{17}$ and $\tilde{r}_3(\pi)=\overline{7}$. For $p=1$, we see that there is a part $\overline{17}$, which is congruent to $7$ modulo $10$, in the $3$-band $\{\tilde{r}_{p+1}(\pi)\}_3=\{\overline{25}_1,\overline{20}_2,\overline{17}_3\}$, and so $\tilde{\tilde{r}}_{p+1}(\pi)=\overline{17}$.  Furthermore, there is no non-degenerate part of $\pi$ less than $\overline{17}$.   Hence, $\pi\in\mathcal{\overline{{C}}}_3(5,7;10,4,3|3,1).$

We consider the overpartition $\pi$ defined in \eqref{new-example-new-1}. Let $\overline{\pi}$ be the overpartition obtained  from $\pi$ by subtracting $3$ from the  non-degenerate 2-part $\pi_{10}=10$. The  Gordon marking of $\overline{\pi}$ is given as follows.
\[G(\overline{\pi})=(\overline{50}_1,\overline{47}_2,40_3,40_1,\overline{30}_2,\overline{25}_1,\overline{20}_3,20_2,\overline{15}_1,\overline{7}_2,\overline{5}_1).\]
There are two parts marked with $3$ in $G(\overline{\pi})$, which are  $\tilde{g}_1(\overline{\pi})=40$ and $\tilde{g}_2(\overline{\pi})=\overline{20}$. So, $\pi\in\mathcal{\overline{{C}}}_{10}(5,7;10,4,3|2,2).$

For the  overpartition $\pi$   defined in \eqref{new-example-11},
we know that $\pi_8=20$ is the smallest non-degenerate part of $\pi$, and so $\pi_m=\pi_8=20$. Let $\overline{\pi}$ be the overpartition obtained  from $\pi$ by subtracting $3$ from the smallest non-degenerate part $\pi_8=20$. The  Gordon marking of $\overline{\pi}$ is illustrated below.
\begin{equation}\label{example-11111}
G(\overline{\pi})=(\overline{50}_2,\overline{47}_1,{ {40}_3},{ {40}_2,\overline{30}_1},{ \overline{25}_3},{ \overline{20}_2,\overline{17}_1},{ \overline{10}_3},{ {10}_2,\overline{7}_1}).
\end{equation}
There are three parts marked with $3$  in  $G(\overline{\pi})$, which are $\tilde{g}_1(\overline{\pi})=40$, $\tilde{g}_2(\overline{\pi})=\overline{25}$ and  $\tilde{g}_3(\overline{\pi})=\overline{10}$. It is easy to  check that  $\tilde{g}_1(\overline{\pi})=40>\pi_m+\eta=30>\tilde{g}_2(\overline{\pi})=\overline{25}$. Hence,
$\pi\in\mathcal{\overline{{C}}}_{10}(5,7;10,4,3|3,1).$

We can express the $(k-1)$-addition introduced by Kim \cite{Kim-2018} in the context of overpartitions as follows.

\begin{defi}[The $(k-1)$-addition]\label{defi-addition-1}
For $0\leq p\leq N,$ let $\pi$ be an overpartition in  $\overline{\mathcal{C}}_\lambda(\alpha_2,\ldots,\alpha_\lambda;\eta,k,r|N,p)$ and let $\tilde{r}_1(\pi)>\cdots>\tilde{r}_N(\pi)$ be the $(k-1)$-marked parts in $RG(\pi)$.  The $(k-1)$-addition  $A_{p\eta+\alpha_1}\colon\pi\rightarrow\tau$ is defined as follows{\rm{:}} First apply the forward move $\phi_p$ to $\pi,$ and then add $\alpha_1$ to $\tilde{\tilde{r}}_{p+1}(\pi)$ to generate a non-overlined part divisible by $\eta$. Here{\rm{,}} we assume that $\tilde{\tilde{r}}_{N+1}(\pi)=\overline{\alpha_\lambda}$. Rearrange the parts in non-increasing order  to obtain the overpartition $\tau$.
\end{defi}

For example, take the overpartition $\pi$ in $\mathcal{\overline{{C}}}_3(5,7;10,4,3|3,1)$, whose reverse Gordon marking is given in \eqref{new-example-21}. Note that $p=1$, we first change the part $\tilde{r}_{1}(\pi)=30$  to $\tilde{r}_{1}(\pi)+\eta=40$ and then add $\alpha_1=3$ to the part $\tilde{\tilde{r}}_{2}(\pi)=\overline{17}$ to get  $20$. So we obtain
\[\tau=A_{1\cdot10+3}(\pi)=(\overline{50},\overline{47},40,40,\overline{30},\overline{25},\overline{20},
20,\overline{10},10,\overline{7}),\]
which is the overpartition defined in \eqref{new-example-11}. Clearly, $f_\tau(0,10]=f_\pi(0,10]=2,$ $|\tau|=|\pi|+p\eta+\alpha_1=|\pi|+13$ and
$\tau\in\mathcal{\overline{{C}}}_{10}(5,7;10,4,3|3,1).$

In light of Proposition 5.1 and Proposition 5.3 in \cite{Kim-2018}, we deduce that the $(k-1)$-addition is a map from $\overline{\mathcal{C}}_\lambda(\alpha_2,\ldots,\alpha_\lambda;\eta,k,r|N,p)$ to $\overline{\mathcal{C}}_\eta(\alpha_2,\ldots,\alpha_\lambda;\eta,k,r|N,p)$.

\begin{lem}\label{dilation-1}
For  $0\leq p\leq N,$  let  $\pi$ be an overpartition in $\overline{\mathcal{C}}_\lambda(\alpha_2,\ldots,\alpha_\lambda;\eta,k,r|N,p)$ and let $\tau=A_{p\eta+\alpha_1}(\pi)$.   Then $\tau$ is an overpartition in $\overline{\mathcal{C}}_\eta(\alpha_2,\ldots,\alpha_\lambda;\eta,k,r|N,p)$ such that $|\tau|=|\pi|+p\eta+\alpha_1.$
 \end{lem}

In the context of overpartitions, the inverse map of the $(k-1)$-addition (i.e., the $(k-1)$-subtraction) can be stated  as follows.

\begin{defi}[The $(k-1)$-subtraction]\label{defi-subtraction-1}
For $0\leq p\leq N$, let $\tau$ be an overpartition in  $\overline{\mathcal{C}}_\eta(\alpha_2,\ldots,\alpha_\lambda;\eta,k,r|N,p)$. The $(k-1)$-subtraction ${S}_{p\eta+\alpha_1}\colon \tau\rightarrow\pi$ is defined as follows{\rm{:}}
First subtract $\alpha_1$ from the smallest non-degenerate part of $\tau$ to get $\overline{\tau}$, and then apply the backward move $\psi_p$ to $\overline{\tau}$ to obtain $\pi$.
\end{defi}

For example, let $\tau$ be the overpartition which agrees with the one in \eqref{new-example-11}. We know that $\tau\in\mathcal{\overline{{C}}}_{10}(5,7;10,4,3|3,1).$ We first subtract $3$ from the smallest non-degenerate part $\tau_{8}=20$ to obtain $\overline{\tau}$, whose Gordon marking is given in \eqref{example-11111}. Then we apply the backward move $\psi_1$ to $\overline{\tau}$ to obtain $\pi$, namely, replace  $\tilde{g}_1(\overline{\tau})=40$ in $\overline{\tau}$ by $30$ in $\pi$. So we get
\[\pi=(\overline{50},\overline{47},40,\overline{30},30,\overline{25},\overline{20},\overline{17},\overline{10},10,\overline{7}),\]
whose revere Gordon marking is given in \eqref{new-example-21}. From the proceeding example \eqref{new-example-21}, we see that $\pi\in\mathcal{\overline{{C}}}_{3}(5,7;10,4,3|3,1).$

By virtue of Proposition 5.2 and Proposition 5.3 in \cite{Kim-2018}, we derive that the $(k-1)$-subtraction is a map from $\overline{\mathcal{C}}_\eta(\alpha_2,\ldots,\alpha_\lambda;\eta,k,r|N,p)$ to $\overline{\mathcal{C}}_\lambda(\alpha_2,\ldots,\alpha_\lambda;\eta,k,r|N,p)$.

\begin{lem}\label{2etatau}
For  $0\leq p\leq N,$  let $\tau$ be an overpartition in $\overline{\mathcal{C}}_\eta(\alpha_2,\ldots,\alpha_\lambda;\eta,k,r|N,p)$ and let $\pi=S_{p\eta+\alpha_1}(\pi)$.   Then $\pi$ is an overpartition in $\overline{\mathcal{C}}_\lambda(\alpha_2,\ldots,\alpha_\lambda;\eta,k,r|N,p)$ such that $|\pi|=|\tau|-p\eta-\alpha_1.$
 \end{lem}

 Combining Lemma \ref{dilation-1} and Lemma \ref{2etatau},  we derive the following result.

\begin{thm}\label{Dilation-reduction}
For $0\leq p\leq N$, the $(k-1)$-addition $A_{p\eta+\alpha_1}$   is a bijection between  $\overline{\mathcal{C}}_\lambda(\alpha_2,\ldots,\alpha_\lambda;
\eta,k,r|N,p)$  and  $\overline{\mathcal{C}}_\eta(\alpha_2,\ldots,\alpha_\lambda;
\eta,k,r|N,p)$. Moreover, assume that $\pi$ is an overpartition in $\overline{\mathcal{C}}_\lambda(\alpha_2,\ldots,\alpha_\lambda;\eta,k,r|N,p),$ let $\tau=A_{p\eta+\alpha_1}(\pi).$ Then, we have $f_\tau(0,\eta]=f_\pi(0,\eta]$ and $|\tau|=|\pi|+p\eta+\alpha_1.$

\end{thm}

We conclude this subsection with the following theorem, which shows that the $(k-1)$-addition operation can be applied to an overpartition in $ \overline{\mathcal{C}}_\lambda(\alpha_2,\ldots,\alpha_\lambda;\eta,k,r|N,p)$ successively.   It is worth mentioning that the following theorem can be obtained from Proposition 5.6 in \cite{Kim-2018}. Here, we provide a detailed proof for completeness.

\begin{thm}\label{tauhh}
 For $0\leq p<N$, let $\pi$ be an overpartition in $ \overline{\mathcal{C}}_\lambda(\alpha_2,\ldots,\alpha_\lambda;\eta,k,r|N,p)$ and let $\tau=A_{p\eta+\alpha_1}(\pi)$.  Assume that there are $N'$ parts marked with $k-1$  in $RG(\tau)$. Then for $p'<N'$,
 $\tau$ is also an overpartition in $\overline{\mathcal{C}}_\lambda(\alpha_2,\ldots,\alpha_\lambda;\eta,k,r|N',p')$ if and only if $p'>p$.

\end{thm}

\pf Since $\pi\in \overline{\mathcal{C}}_\lambda(\alpha_2,\ldots,\alpha_\lambda;\eta,k,r|N,p)$,  there are $N$ parts marked with $k-1$ in $RG(\pi)$, denoted  $\tilde{r}_1(\pi)>\cdots>\tilde{r}_N(\pi)$. Furthermore, there exists a part $\equiv \alpha_\lambda \pmod{\eta}$ in the $(k-1)$-band  $\{\tilde{r}_{p+1}(\pi)\}_{k-1}$,  denoted  $\tilde{\tilde{r}}_{p+1}(\pi)$, and there is no non-degenerate part of $\pi$ less than $\tilde{\tilde{r}}_{p+1}(\pi)$. Appealing to Theorem  \ref{Dilation-reduction}, we know that $\tau\in\overline{\mathcal{C}}_\eta(\alpha_2,\ldots,
\alpha_\lambda;\eta,k,r|N,p)$ and $\tilde{\tilde{r}}_{p+1}(\pi)+\alpha_1$ is the smallest non-degenerate part of $\tau$.

We first show that if $N'>p'>p$, then $\tau$ is in $\overline{\mathcal{C}}_\lambda(\alpha_2,\ldots,\alpha_\lambda;\eta,k,r|N',p')$. Let $\tilde{r}_1(\tau)>\cdots>\tilde{r}_{N'}(\tau)$ be the $(k-1)$-marked parts in $RG(\tau)$.
We are required to prove that there exists a part $\equiv \alpha_\lambda \pmod{\eta}$ in the $(k-1)$-band $\{\tilde{r}_{p'+1}(\tau)\}_{k-1}$ of $\tau$, denoted $\tilde{\tilde{r}}_{p'+1}(\tau)$,
and there is no non-degenerate part of $\tau$ less than $\tilde{\tilde{r}}_{p'+1}(\tau)$. Note that $\tilde{\tilde{r}}_{p+1}(\pi)+\alpha_1$ is the smallest non-degenerate part of $\tau$, it suffices to prove that the parts in the $(k-1)$-band $\{\tilde{r}_{p'+1}(\tau)\}_{k-1}$ are less than $\tilde{\tilde{r}}_{p+1}(\pi)+\alpha_1.$ To do this, under the assumption that $p'>p$, it is enough to verify that
\begin{equation}\label{proof-ss-1}
\tilde{r}_{p+1}(\tau)\leq \tilde{\tilde{r}}_{p+1}(\pi)+\alpha_1.
\end{equation}
Suppose to the contrary that $\tilde{r}_{p+1}(\tau)>\tilde{\tilde{r}}_{p+1}(\pi)+\alpha_1.$ By the definition of reverse Gordon marking, we see that there is a $(k-1)$-band of $\tau$ in $(\tilde{\tilde{r}}_{p+1}(\pi)+\alpha_1,\tilde{r}_p(\tau))$.

Assume that $\pi^{(1)}=\phi_p(\pi)$. By virtue of (3) and (4) in Proposition 4.3 in \cite{he-ji-zhao}, we know  that there are $N$ parts marked with $k-1$ in $G(\pi^{(1)})$, denoted $\tilde{g}_1(\pi^{(1)})>\cdots>\tilde{g}_N(\pi^{(1)})$, and $\tilde{g}_p(\pi^{(1)})= \tilde{r}_p(\pi)+\eta$. By Proposition \ref{sequence-length}, we find that there are also $N$ parts marked with $k-1$ in $RG(\pi^{(1)})$, denoted  $\tilde{r}_1(\pi^{(1)})>\cdots>\tilde{r}_N(\pi^{(1)})$, and $\tilde{g}_p(\pi^{(1)})-\eta\leq \tilde{r}_p(\pi^{(1)})\leq \tilde{g}_p(\pi^{(1)})$. So, we get $\tilde{r}_p(\pi)\leq\tilde{r}_p(\pi^{(1)})\leq \tilde{r}_p(\pi)+\eta$. Recall that $\tau$ is obtained by replacing $\tilde{\tilde{r}}_{p+1}(\pi)$ in $\pi^{(1)}$ by  $\tilde{\tilde{r}}_{p+1}(\pi)+\alpha_1$, so $\tilde{r}_p(\tau)=\tilde{r}_p(\pi^{(1)})$.
It yields
\begin{equation}\label{aaaaaa-11}
\tilde{r}_p(\pi)\leq\tilde{r}_p(\tau)\leq \tilde{r}_p(\pi)+\eta.
 \end{equation}
 Under the condition that there is a $(k-1)$-band of $\tau$ in $(\tilde{\tilde{r}}_{p+1}(\pi)+\alpha_1,\tilde{r}_p(\tau))$, we deduce that there is a $(k-1)$-band of $\tau$ in $(\tilde{\tilde{r}}_{p+1}(\pi)+\alpha_1,\tilde{r}_p(\pi)+\eta)$.

 Let $\{\tau_{i+l}\}_{0\leq l\leq k-2}$   be a $(k-1)$-band of $\tau$ in $(\tilde{\tilde{r}}_{p+1}(\pi)+\alpha_1,
\tilde{r}_{p}(\pi)+\eta)$.  Recall that  $\tilde{r}_{p}(\pi)$ in $\pi$ is changed to  $\tilde{r}_{p}(\pi)+\eta$ in $\tau$, so  $\tau_i<\tilde{r}_{p}(\pi)$. This implies that  $\{\tau_{i+l}\}_{0\leq l\leq k-2}$  is a $(k-1)$-band of $\tau$ in $(\tilde{\tilde{r}}_{p+1}(\pi)+\alpha_1,\tilde{r}_{p}(\pi))$. From the construction of $\tau$, we find that the parts in the $(k-1)$-band $\{\tau_{i+l}\}_{0\leq l\leq k-2}$ of $\tau$ are also parts of $\pi$.  It yields that $\{\tau_{i+l}\}_{0\leq l\leq k-2}$  is also a $(k-1)$-band of $\pi$ in $(\tilde{\tilde{r}}_{p+1}(\pi)+\alpha_1,\tilde{r}_{p}(\pi))$, which leads to a contradiction.  Hence, we arrive at \eqref{proof-ss-1}. This completes the proof of the sufficiency.

Conversely, assume that $\tau$ is also an overpartition in $\overline{\mathcal{C}}_\lambda(\alpha_2,\ldots,\alpha_\lambda;\eta,k,r|N',p')$, and for $p'<N'$,  we intend to show that $p'>p$.  Suppose to the contrary that $p'\leq p$.

Given that $\tau\in\overline{\mathcal{C}}_\lambda(\alpha_2,\ldots,\alpha_\lambda;\eta,k,r|N',p')$, then there exists a part $\equiv \alpha_\lambda \pmod{\eta}$ in the $(k-1)$-band $\{\tilde{r}_{p'+1}(\tau)\}_{k-1}$ of $\tau$, denoted $\tilde{\tilde{r}}_{p'+1}(\tau)$, and there is no non-degenerate part of $\tau$ less than $\tilde{\tilde{r}}_{p'+1}(\tau)$. Recall that  $\tilde{\tilde{r}}_{p+1}(\pi)+\alpha_1$ is a non-degenerate part of $\tau$, we have $\tilde{\tilde{r}}_{p'+1}(\tau)\leq\tilde{\tilde{r}}_{p+1}(\pi)+\alpha_1$.
Under the condition that $\tilde{\tilde{r}}_{p'+1}(\tau)$ and $\tilde{\tilde{r}}_{p+1}(\pi)$ are both congruent to $\alpha_\lambda$ modulo $\eta$, we deduce that $\tilde{\tilde{r}}_{p'+1}(\tau)\leq\tilde{\tilde{r}}_{p+1}(\pi)$.   Furthermore, from the construction of $\tau$, we see that  $\tilde{\tilde{r}}_{p+1}(\pi)$  does not occur in $\tau$. It yields $\tilde{\tilde{r}}_{p'+1}(\tau)\leq\tilde{\tilde{r}}_{p+1}(\pi)-\eta$. Since $\tilde{\tilde{r}}_{p+1}(\pi)-\eta<\tilde{r}_{p+1}(\pi)$, we get
\begin{equation}\label{contradiction-new}
\tilde{\tilde{r}}_{p'+1}(\tau)<\tilde{r}_{p+1}(\pi).
\end{equation}

On the other hand, appealing to \eqref{aaaaaa-11}, we get $\tilde{r}_p(\tau)\geq\tilde{r}_p(\pi)$.
 Assume that $\tilde{r}_{p+1}(\pi)$ is the $m$-th part of $\pi$, that is, $\pi_m=\tilde{r}_{p+1}(\pi)$. Since $\{\pi_{m-l}\}_{0\leq l\leq k-2}$ is the $(k-1)$-band of $\pi$ induced by $\tilde{r}_{p+1}(\pi)$, we could assume that $\pi_{m-t}=\tilde{\tilde{r}}_{p+1}(\pi)$, where $0\leq t\leq k-2$. It is easy to check that
$\{\pi_{m-k+2},\ldots,\pi_{m-t}+\alpha_1,\ldots,\pi_m\}$
is a $(k-1)$-band of $\tau$. It follows that $\tilde{r}_{p+1}(\tau)\geq\tilde{r}_{p+1}(\pi)$. Under the assumption that $p'\leq p$, we get $\tilde{\tilde{r}}_{p'+1}(\tau)\geq\tilde{r}_{p'+1}(\tau)\geq\tilde{r}_{p+1}(\tau)\geq\tilde{r}_{p+1}(\pi)$, which is in contradiction to \eqref{contradiction-new}. Thus, we have shown  $p'>p$.
This completes the proof.       \qed

  \subsection{The $(k-1)$-insertion and the $(k-1)$-separation}

 The main objective of this subsection is    to
 recall  the $(k-1)$-insertion with $a=\alpha_1$ defined in \cite{he-ji-zhao}, which can be viewed as an  overpartition generalization of the $(k-1)$-insertion  introduced by Kim \cite{Kim-2018}. It is worth mentioning that we could not apply the $(k-1)$-insertion  introduced by Kim   \cite{Kim-2018} directly in the second step owing to the presence of certain parts, such as $t\eta<\overline{t\eta}<\overline{t\eta+\alpha_1}$ in the overpartition belonging to  $\mathcal{\overline{{B}}}_1(\alpha_1,\ldots,\alpha_\lambda;\eta,k,r)$.

  The $(k-1)$-insertion with $a=\alpha_1$ can be used to merge  the remaining parts of $\delta^{(1)}$ (which are $\equiv \alpha_1\pmod{\eta}$) and the   overpartition in $\mathcal{\overline{{C}}}_1(\alpha_2,\ldots,\alpha_\lambda;\eta,k,r)$ to generate certain overlined parts $\equiv \alpha_1\pmod{\eta}$, which results in   an overpartition in $\mathcal{\overline{{B}}}_1(\alpha_1,\ldots,\alpha_\lambda;\eta,k,r)$.  We first recall two subsets of $\mathcal{\overline{{B}}}_1(\alpha_1,\ldots,\alpha_\lambda;\eta,k,r)$ introduced in \cite{he-ji-zhao}.

 \begin{itemize}
 \item    For  $s\geq N\geq 0$, let $\overline{\mathcal{B}}^{\, {\alpha_1}}_<(\alpha_1,\ldots,\alpha_\lambda;\eta,k,r|N,s)$ denote the set of   overpartitions $\tau$ in $\overline{\mathcal{B}}_1(\alpha_1,\ldots,\alpha_\lambda;\eta,k,r)$  satisfying
   \begin{itemize}
 \item[{\rm (1)}] There are $N$ parts marked with $k-1$ in   $RG(\tau)$, denoted $ \tilde{r}_1(\tau)>\cdots > \tilde{r}_N(\tau)${\rm{;}}

  \item[{\rm (2)}]  Assume that  $p$ is the smallest integer satisfying  $ \tilde{r}_{p+1}(\tau)+\eta \leq \overline{(s-p)\eta+{\alpha_1}}$ with the convention that $\tilde{r}_{N+1}(\tau)=-\infty$. Then the largest overlined part $\equiv {\alpha_1}\pmod\eta$ in $\tau$ is less than $\overline{(s-p)\eta+{\alpha_1}}${\rm{;}}

 \item[{\rm (3)}]   If  $f_\tau(0,\eta]=r-1$ and $s=N\geq 1$, then  $\tilde{r}_{N}(\tau)\leq\eta${\rm{;}}

      \item[{\rm (4)}]  If   $s=N=0$, then $f_\tau(0,\eta]<r-1$.

  \end{itemize}

  \item For $s\geq N\geq 0$, let $\overline{\mathcal{B}}^{\, {\alpha_1}}_{=}(\alpha_1,\ldots,\alpha_\lambda;\eta,k,r|N,s)$ denote the set of overpartitions $\sigma$ in $\overline{\mathcal{B}}_1(\alpha_1,\ldots,\alpha_\lambda;\eta,k,r)$    subject to the following conditions:
     \begin{itemize}
  \item[{\rm (1)}] There exists  an  overlined part $\equiv \alpha_1\pmod\eta$ in $\sigma$, and assume that the  largest overlined part $\equiv \alpha_1\pmod\eta$  in  $\sigma$ is $\overline{t\eta+\alpha_1}$;

    \item[{\rm (2)}] Let ${\hat{\sigma}}$ be the  overpartition obtained  by removing
  $\overline{t\eta+\alpha_1}$  from $\sigma$.  Then there are $N$ parts marked with $k-1$ in  $G({\hat{\sigma}})$,  denoted $ \tilde{g}_1({\hat{\sigma}})>\cdots > \tilde{g}_N({\hat{\sigma}})$;

    \item[{\rm (3)}] Assume that  $p$  is the smallest integer such that  $ \tilde{g}_{p+1}({\hat{\sigma}})< \overline{t\eta+\alpha_1}$ with the convention that $\tilde{g}_{N+1}({\hat{\sigma}})=-\infty$.  Then $s=p+t$.

      \end{itemize}

\end{itemize}

For example, let $s=6$ and $\alpha_1=1$ and let $\tau$ be the overpartition in $\overline{\mathcal{B}}_1(1,5,9;10,5,4)$ with the reverse Gordon marking
  \begin{equation}\label{example-insertion-over-1}
\begin{split}
&RG(\tau)=(\overline{85}_1,{ {80}_2,{80}_3,\overline{75}_1},{ {70}_4}
,\overline{69}_2,
{ \overline{60}_1,{60}_3,\overline{59}_2}, { \overline{55}_4},{ \overline{50}_1,\overline{49}_2,\overline{45}_3,\overline{40}_1,{40}_4},\\[5pt]
&\ \ \ \ \ \ \ \ \ \ \ \ \ \  \overline{39}_2,\overline{35}_3,{{30}_1, \overline{29}_2,\overline{25}_3,} { \overline{21}_4},
{\overline{19}_1, \overline{15}_2,\overline{11}_3,{10}_4,\overline{9}_1},
{ \overline{5}_2}).
\end{split}
\end{equation}
There are five  $4$-marked parts  $\tilde{r}_1(\tau)={70}$, $\tilde{r}_2(\tau)=\overline{55}$, $\tilde{r}_3(\tau)={40}$, $\tilde{r}_4(\tau)=\overline{21}$ and $\tilde{r}_5(\tau)=10$  in  $RG(\tau)$.   Then $p=3$ is the smallest integer such that $ \overline{31}= \tilde{r}_{p+1}(\tau)+\eta\leq \overline{(s-p)\eta+\alpha_1}=\overline{31}$. Meanwhile,  the largest overlined part $\equiv {1}\pmod{10}$ in $\tau$ is $\overline{21}$, which is less than $\overline{(s-p)\eta+\alpha_1}=\overline{31}$.  So $\tau$ is an overpartition in $\overline{\mathcal{B}}^{\,1}_<(1,5,9;10,5,4|5,6)$.

For another example, let $s=6$ and $\alpha_1=1$ and let
  \begin{equation}\label{example-insertion-over-2}
\begin{split}
&\sigma=(\overline{85},{{80}}
,{ {80},{80},\overline{75}},\overline{69},
 {\overline{65}},{ \overline{60},{60},\overline{59}},{\overline{50}},{ {50},\overline{49},\overline{45}},{ \overline{40}},\\[5pt]
&\ \ \ \ \ \ \
\overline{39},\overline{35},\overline{31},{30},{ \overline{29},\overline{25},\overline{21},}
{ \overline{19},\overline{15},\overline{11}},10,
{\overline{9}},\overline{5})
\end{split}
\end{equation}
  be an overpartition in $\overline{\mathcal{B}}_1(1,5,9;10,5,4)$.  The largest overlined part $\equiv {1}\pmod{10}$ in $\sigma$ is  $\overline{31}$, and so $t=3$. Removing  $\overline{31}$  from $\sigma$, we get ${\hat{\sigma}}$ with the Gordon marking  \begin{equation}\label{example-insertion-over-3}
\begin{split}
&G({\hat{\sigma}})=(\overline{85}_2,{ {80}_4}
,{ {80}_3,{80}_1,\overline{75}_2},\overline{69}_1,
 { \overline{65}_4},{ \overline{60}_3,{60}_2,\overline{59}_1},{ \overline{50}_3},{ {50}_4,\overline{49}_2,\overline{45}_1,\overline{40}_3},\\[5pt]
&\ \ \ \ \ \ \ \ \ \ \ \
\overline{39}_2,\overline{35}_1,{{30}_4,\overline{29}_2,{ \overline{25}_1},} { \overline{21}_3},{\overline{19}_2},
{ \overline{15}_1,{ \overline{11}_4},{10}_3},
{\overline{9}_2},\overline{5}_1).
\end{split}
\end{equation}
There are five  $4$-marked parts  $\tilde{g}_1({\hat{\sigma}})={80}$, $\tilde{g}_2({\hat{\sigma}})=\overline{65}$, $\tilde{g}_3({\hat{\sigma}})=50$, $\tilde{g}_4({\hat{\sigma}})=30$ and $\tilde{g}_5({\hat{\sigma}})=\overline{11}$  in  $G({\hat{\sigma}})$ and $p=3$ is the smallest integer such that
$30=\tilde{g}_{p+1}({\hat{\sigma}})<\overline{31}$. Indeed, $p+t=s$ holds.   Thus, we conclude that $\sigma$ is an overpartition in $\overline{\mathcal{B}}^{\,1}_=(1,5,9;10,5,4|5,6)$.

The following is the definition of the $(k-1)$-insertion.

\begin{defi}{\rm\cite[Definition 4.9]{he-ji-zhao}}\label{the insertion} For  $s\geq N\geq 0,$  assume that $\tau$ is an overpartition in $\overline{\mathcal{{B}}}^{\, \alpha_1}_<(\alpha_1,\ldots,\alpha_\lambda;\eta,k, r|N,s)$ with $N$ parts marked with $k-1$  in $RG(\tau),$ denoted  $\tilde{r}_1(\tau)>\cdots>\tilde{r}_N(\tau)$. Let  $p$ be the smallest integer such that   $0\leq p\leq N$ and $\overline{(s-p)\eta+\alpha_1}\geq\tilde{r}_{p+1}(\tau)+\eta$. The $(k-1)$-insertion $I^{\alpha_1}_s\colon \tau \rightarrow \sigma$ is defined  as follows{\rm :}   first apply the forward move $\phi_p$ to $\tau$ to get ${{\tau}'}$,  then insert $\overline{(s-p)\eta+\alpha_1}$ into $\tau'$  as an overlined part of $\sigma$.
 \end{defi}

For example, let  $\tau$ be the overpartition in   $\overline{\mathcal{B}}^{\,1}_<(1,5,9;10,5,4|5,6)$ whose reverse Gordon marking is  given in \eqref{example-insertion-over-1}. In this case,  $p=3$   is the smallest integer such that $\overline{(s-p)\eta+\alpha_1}=\overline{31}\geq \tilde{r}_{p+1}(\tau)+\eta=\overline{31}$.   Applying the forward move  $\phi_3$ to $\tau$, we  get
 \begin{equation*}
\begin{split}
&{\tau'}=(\overline{85}, {{80}}
,{ {80}, {80}, \overline{75}}, \overline{69},
 {\overline{65}}, { \overline{60}, {60}, \overline{59}},{\overline{50}}, { {50},\overline{49},  \overline{45}},\overline{40}, \\[5pt]
&\ \ \ \ \ \ \ \
\overline{39}, \overline{35},   {30, \overline{29}, \overline{25},} { \overline{21}}, \overline{19},
{ \overline{15}, \overline{11}, {10}},
{\overline{9}}, \overline{5}),
\end{split}
\end{equation*}
 whose Gordon marking agrees with the one in \eqref{example-insertion-over-3}. Inserting $\overline{(s-p)\eta+\alpha_1}=\overline{31}$  into  ${{\tau}'}$, we obtain $\sigma=I_s^{\alpha_1}(\tau)$ as in \eqref{example-insertion-over-2}. Clearly,  $|\sigma|=|\tau|+61$.

We find that the $(k-1)$-insertion $I^{\alpha_1}_s$ is a map from  $\overline{\mathcal{{B}}}^{\,\alpha_1}_<(\alpha_1,\ldots,\alpha_\lambda;\eta,k, r|N,s)$ to $\overline{\mathcal{{B}}}^{\,\alpha_1}_{=}(\alpha_1,\ldots,\alpha_\lambda;\eta,k, r|N,s)$.

\begin{lem}{\rm\cite[Lemma 4.11]{he-ji-zhao}}\label{isub-lem-a} For  $s\geq N\geq 0$,
assume that $\tau$ is an overpartition in $\overline{\mathcal{{B}}}^{\,\alpha_1}_{<}(\alpha_1,\ldots,\alpha_\lambda;\eta,k, r|N,s)$ and let $\sigma={I}^{\alpha_1}_{s}(\tau)$. Then  $\sigma$ is an overpartition in $\overline{\mathcal{{B}}}^{\,\alpha_1}_=(\alpha_1,\ldots,\alpha_\lambda;\eta,k, r|N,s)$. Moreover,  $|\sigma|=|\tau|+s\eta+\alpha_1$.
\end{lem}

We now define the $(k-1)$-separation, which is the inverse map of the $(k-1)$-insertion.

\begin{defi}{\rm\cite[Definition 4.12]{he-ji-zhao}}\label{isub} For  $s\geq N\geq 0$,
let $\sigma$ be an overpartition in $\overline{\mathcal{{B}}}^{\,\alpha_1}_{=}(\alpha_1,\ldots,\alpha_\lambda;\eta,k, r|N,s)$ with   the largest overlined part $\equiv \alpha_1\pmod \eta$ being $\overline{t\eta+\alpha_1}$.   The $(k-1)$-separation map ${J}^{ \alpha_1}_{s}\colon \sigma \rightarrow \tau$ is defined as follows{\rm :} First remove  $\overline{t\eta+\alpha_1}$ from $\sigma$ to produce ${\hat{\sigma}}$, and then apply the backward move $\psi_{s-t}$ to ${\hat{\sigma}}$  to obtain $\tau$.
\end{defi}

 The following lemma states that the $(k-1)$-separation is  a map from  $\overline{\mathcal{{B}}}^{\,\alpha_1}_{=}(\alpha_1,\ldots,\alpha_\lambda;\eta,\break$
 $k, r|N,s)$ to $\overline{\mathcal{{B}}}^{\,\alpha_1}_<(\alpha_1,\ldots,\alpha_\lambda;\eta,k, r|N,s)$.

\begin{lem}{\rm\cite[Lemma 4.13]{he-ji-zhao}}\label{isub-lem} For  $s\geq N\geq 0$,
assume that $\sigma$ is an overpartition in $\overline{\mathcal{{B}}}^{\,\alpha_1}_{=}(\alpha_1,\ldots,\alpha_\lambda;\eta,k, r|N,s)$ and let $\tau={J}^{\alpha_1}_{s}(\sigma)$. Then  $\tau$ is an overpartition in $\overline{\mathcal{{B}}}^{\,\alpha_1}_<(\alpha_1,\ldots,\alpha_\lambda;\eta,k, r|N,s)$. Moreover,  $|\tau|=|\sigma|-s\eta-\alpha_1$.
\end{lem}

Based on Lemma \ref{isub-lem-a}  and Lemma \ref{isub-lem}, we reach the following consequence.

\begin{thm}{\rm \cite[Theorem 4.10]{he-ji-zhao}}\label{deltagammathmbb} For  $s\geq N\geq 0$,  the $(k-1)$-insertion  $ I^{\alpha_1}_{s}$   is a bijection between   $\mathcal{\overline{B}}^{\, \alpha_1}_<(\alpha_1,\ldots,\alpha_\lambda;\eta,k, r|N,s)$  and   $\mathcal{\overline{B}}^{\, \alpha_1}_=(\alpha_1,\ldots,\alpha_\lambda;\eta,k, r|N,s)$. Moreover, for $\tau \in \overline{\mathcal{{B}}}^{\, \alpha_1}_<(\alpha_1,\ldots,\alpha_\lambda;\eta,k, r|N,s)$, let $\sigma=I^{\alpha_1}_{s}(\tau)$, then we have $|\sigma|=|\tau|+s\eta+\alpha_1$.
 \end{thm}

The following theorem gives a   criterion   to determine whether an overpartition in $\overline{\mathcal{{B}}}^{\,\alpha_1}_{=}(\alpha_1,\ldots,\alpha_\lambda;\eta,k, r|N,s)$ is also an overpartition in $\overline{\mathcal{{B}}}^{\,\alpha_1}_{<}(\alpha_1,\ldots,\alpha_\lambda;\eta,k, r|N',s')$, which involves  the successive application of the $(k-1)$-insertion operation.

\begin{thm}{\rm \cite[Theorem 4.14]{he-ji-zhao}}\label{ssins}
For $s\geq N\geq 0$, let $\sigma$ be an overpartition in $\overline{\mathcal{{B}}}^{\,\alpha_1}_{=}(\alpha_1,\ldots,\alpha_\lambda;\eta,k, r|N,s)$. Assume that there are $N'$ parts marked with $k-1$ in  $RG(\sigma)$. Then,  $\sigma$ is also an overpartition in $\overline{\mathcal{{B}}}^{\,\alpha_1}_<(\alpha_1,\ldots,\alpha_\lambda;\eta,k, r|N',s')$ if and only if $s'>s$.

\end{thm}

\subsection{Proof of Theorem \ref{lambdathm}}

In this subsection, we will give a proof of Theorem \ref{lambdathm} by successively applying the   $(k-1)$-addition and the $(k-1)$-insertion.   Before approaching the proof of Theorem \ref{lambdathm}, we first prove
the following theorem, which tells us that   the $(k-1)$-insertion can be applied to the resulting overpartition obtained by applying the $(k-1)$-addition.

 \begin{thm}\label{sskr1}
Let $\tau$ be an overpartition in $\overline{\mathcal{C}}_\eta(\alpha_2,\ldots,\alpha_\lambda;\eta,k,r|N,N)$. Assume that there are $N'$ parts marked with $k-1$ in $RG(\tau)$. Then $\tau\in\overline{\mathcal{B}}^{\,\alpha_1}_<(\alpha_1,\ldots,\alpha_\lambda;\eta,k,r|N',s)$ if and only if $s>N$.
\end{thm}

 \pf  Let $\tilde{r}_{1}(\tau)>\cdots>\tilde{r}_{N'}(\tau)$ be the $(k-1)$-marked parts in $RG(\tau)$. By definition, we see that  $N'=N$ or $N+1$.
 Let $\overline{\tau}$ be the  overpartition obtained from $\tau$ by subtracting   $\alpha_1$ from the   non-degenerate $(r-1)$-part $\eta$ of $\tau$ and  let $\tilde{r}_{1}(\overline{\tau})>\cdots>\tilde{r}_{N}(\overline{\tau})$ be the $(k-1)$-marked parts in $RG(\overline{\tau})$.

  We proceed to show that if $s>N$, then $\tau\in\overline{\mathcal{B}}^{\,\alpha_1}_<(\alpha_1,\ldots,\alpha_\lambda;\eta,k,r|N',s)$. Recall that $N'=N$ or $N+1$, then we have $s\geq N+1\geq N'$.
  Assume that $p$ is the smallest integer such that $0\leq p\leq N'$ and $\tilde{r}_{p+1}(\tau)+\eta\leq \overline{(N'-p)\eta+\alpha_1}$.  Note that there are no overlined parts $ \equiv\alpha_1\pmod\eta$ in $\tau$, so the overlined parts $\equiv\alpha_1\pmod\eta$ in $\tau$ are less than $\overline{(N'-p)\eta+\alpha_1}$.

 Since $\tau\in\overline{\mathcal{C}}_\eta(\alpha_2,\ldots,\alpha_\lambda;\eta,k,r|N,N)$, we have $f_\tau(0,\eta]=r-1$. If $s=N'$,   then we have $N'=N+1\geq 1$. Note that there are $N$ parts marked with $k-1$ in $RG(\overline{\tau})$, we see that   $\tilde{r}_{N'}(\tau)=\tilde{r}_{N+1}(\tau)\leq \eta$. Hence we arrive at $\tau\in\overline{\mathcal{B}}^{\,\alpha_1}_<(\alpha_1,\ldots,\alpha_\lambda;\eta,k,r|N',s)$. This completes the proof of the sufficiency.

  Conversely, assume that $\tau\in\overline{\mathcal{B}}^{\,\alpha_1}_<(\alpha_1,\ldots,\alpha_\lambda;\eta,k,r|N',s)$, where $s\geq N'$. We intend to prove that $s > N$.  Suppose to the contrary that $s=N'=N$. In this case, $\tau\in\overline{\mathcal{B}}^{\,\alpha_1}_<(\alpha_1,\ldots,\alpha_\lambda;\eta,k,r|N,N)$.
  It follows from $\tau\in\overline{\mathcal{C}}_\eta(\alpha_2,\ldots,\alpha_\lambda;\eta,k,r|N,N)$ that  $f_\tau(0,\eta]=r-1$. By the definition of $\overline{\mathcal{B}}^{\,\alpha_1}_<(\alpha_1,\ldots,\alpha_\lambda;\eta,k,r|N,N)$, we get $s=N\geq 1$ and $\tilde{r}_N({\tau})\leq\eta$. But, again by  $\tau\in\overline{\mathcal{C}}_\eta(\alpha_2,\ldots,\alpha_\lambda;\eta,k,r|N,N)$, we see that there are no $(k-1)$-bands of $\tau$ in $(0,\overline{\eta+\alpha_\lambda})$, which implies that $\tilde{r}_{N'}(\tau)\geq \eta$. Moreover, we have $\tilde{r}_{N}(\tau)=\tilde{r}_{N}(\overline{\tau})>\eta$,  which leads  to a contradiction.
  Thus, we have shown $s > N$. This completes the proof.   \qed

We are now in a position to prove the main result of this paper.

\noindent{\bf Proof of Theorem \ref{lambdathm}:}  Let $\pi$ be an overpartition in $\mathcal{\overline{B}}_1(\alpha_{2},\ldots,\alpha_{\lambda-1};\eta,k-1,r-1)$ and let $\delta^{(1)}$ be a partition with distinct parts $\equiv \alpha_1 \pmod{\eta}$ and let $\delta^{(\lambda)}$ be a   partition with distinct parts $\equiv \alpha_\lambda \pmod{\eta}$.
We wish to construct an overpartition  $\tau=\Theta(\delta^{(1)},\delta^{(\lambda)},\pi)$  in $\mathcal{\overline{B}}_1(\alpha_{1},\ldots,\alpha_{\lambda};\eta,k,r)$ such that $|\tau|=|\pi|+|\delta^{(1)}|+|\delta^{(\lambda)}|$.

We first insert the parts of $\delta^{(\lambda)}$  as overlined parts  into $\pi$, and denote the resulting overpartition by $\pi^{(0)}$. Evidently, $\pi^{(0)}$ is an  overpartition in $\mathcal{\overline{C}}_1(\alpha_{2},\ldots,
\alpha_{\lambda};\eta,k,r)$ such that
 \begin{equation}\label{thm6.1aaa}
 |\pi^{(0)}|=|\pi|+|\delta^{(\lambda)}|.
 \end{equation}
 Moreover, there is a part $\equiv \alpha_\lambda \pmod{\eta}$ in each  $(k-1)$-band of $\pi^{(0)}$, and $\overline{\alpha_\lambda}$ is a part of $\pi^{(0)}$ when $f_{\pi^{(0)}}(0,\eta]=r-1$. This implies that there is no non-degenerate part of $\pi^{(0)}$. To construct $\tau=\Theta(\delta^{(1)},\delta^{(\lambda)},\pi)$,  we consider  the following two  cases.

\noindent  Case 1: $\delta^{(1)}=\emptyset$. Set $\tau=\pi^{(0)}$. Obviously, $\tau\in\mathcal{\overline{B}}_1(\alpha_{1},\ldots,
\alpha_{\lambda};\eta,k,r)$ and $|\tau|=|\pi|+|\delta^{(1)}|+|\delta^{(\lambda)}|$.

\noindent  Case 2: $\delta^{(1)}\neq\emptyset$. Set $\delta^{(1)}=(q_1\eta+\alpha_1,\ldots,q_m\eta+\alpha_1)$, where $q_1>\cdots>q_m\geq0$. We plan to merge $q_{m}\eta+\alpha_1,\ldots,q_{1}\eta+\alpha_1$ into $\pi^{(0)}$ by performing the $(k-1)$-addition operation and the $(k-1)$-insertion operation successively. The method consists of three steps, and we will indicate the resulting pairs as $({\rm Step}_i(\delta^{(1)}),{\rm Step}_i(\pi^{(0)}))$ after the $i$-th step, where $i=1,2,3$.

 {\bf Step 1:}  We first merge some parts  of $\delta^{(1)}$ from smallest to largest into $\pi^{(0)}$ by  successively applying  the $(k-1)$-addition operation. We denote the intermediate overpartitions by $\pi^{(0)},\, \pi^{(1)}$, and so on. Assume that there are $N(\pi^{(i)})$ parts marked with $k-1$ in $RG(\pi^{(i)})$ for $i\geq 0$. If  $q_m<N(\pi^{(0)})$, note that there is no non-degenerate part of $\pi^{(0)}$,  then we have
\[\pi^{(0)}\in\overline{\mathcal{C}}_\lambda(\alpha_{2},\ldots,\alpha_{\lambda};\eta,k,r|N(\pi^{(0)}),
q_m).\]

Set $b=0$ and repeat the following procedure until $q_{m-b}\geq N(\pi^{(b)})$.

(A)
Apply the $(k-1)$-addition  $A_{q_{m-b}\eta+\alpha_1}$ to $\pi^{(b)}$ to get $\pi^{(b+1)}$, that is,
    \[\pi^{(b+1)}=A_{q_{m-b}\eta+\alpha_1}(\pi^{(b)}).\]
    Since
     \[\pi^{(b)}\in\overline{\mathcal{C}}_\lambda(\alpha_{2},\ldots,\alpha_{\lambda};
\eta,k,r|N(\pi^{(b)}),q_{m-b}),\]
  by Lemma \ref{dilation-1}, we see that
    \begin{equation}\label{ddt}    \pi^{(b+1)}\in\overline{\mathcal{C}}_\eta(\alpha_{2},\ldots,\alpha_{\lambda};\eta,k,r|N(\pi^{(b)}),q_{m-b})
    \end{equation}
    and
    \begin{equation}\label{dda}
    |\pi^{(b+1)}|=|\pi^{(b)}|+q_{m-b}\eta+\alpha_1.
     \end{equation}

(B) Replace $b$ by $b+1$. If $q_{m-b}\geq N(\pi^{(b)})$, then we are done. If $q_{m-b}<N(\pi^{(b)})$, then by  $q_{m-b}>q_{m-b+1}$, we deduce from  Theorem \ref{tauhh}  that
    \[\pi^{(b)}\in\overline{\mathcal{C}}_\lambda(\alpha_{2},\ldots,\alpha_{\lambda};
\eta,k,r|N(\pi^{(b)}),q_{m-b}).\]
Go back to (A).

It is clear that the above procedure terminates after at most $m$ iterations. Assume that the above process terminates with $b=m-j$, that is,  $q_j\geq N(\pi^{(m-j)})$. Set
\[{\rm Step}_1(\pi^{(0)})=\pi^{(m-j)} \quad \text{and} \quad
 {\rm Step}_1(\delta^{(1)})=(q_1\eta+\alpha_1,\ldots,q_j\eta
+\alpha_1).\]
By \eqref{ddt}   and \eqref{dda}, we deduce that
\[{\rm Step}_1(\pi^{(0)})\in\overline{\mathcal{C}}_\eta(\alpha_{2},\ldots,\alpha_{\lambda};\eta,k,r|N(\pi^{(m-j-1)}),q_{j+1})\]
and
\begin{equation}\label{sum-m-t}
|{\rm Step}_1(\pi^{(0)})|=|\pi^{(0)}|+(q_{m}\eta+\alpha_1)+\cdots+(q_{j+1}\eta+\alpha_1).
\end{equation}

{\bf Step 2}: Set  $\sigma={\rm Step}_1(\pi^{(0)})=\pi^{(m-j)}$, and assume that there are  $N(\sigma)$ parts marked with $k-1$ in  $RG(\sigma)$. Recall that $q_j\geq N(\pi^{(m-j)})$, that is  $q_j\geq N(\sigma)$. We consider the following two cases:

Case 1: $q_j=N(\sigma)$ and $\sigma\in\overline{\mathcal{C}}_\lambda(\alpha_{2},\ldots,\alpha_{\lambda};\eta,k,r|q_j,
q_j).$  Then  apply the $(k-1)$-addition $A_{q_{j}\eta+\alpha_1}$ to merge $q_{j}\eta+\alpha_1$ into $\sigma$. In this case, set
\begin{equation}\label{main-step2a}  {\rm Step}_2(\pi^{(0)})=A_{q_{j}\eta+\alpha_1}(\sigma)
\quad  \text{and} \quad
{\rm Step}_2(\delta^{(1)})=(q_1\eta+\alpha_1,\ldots,q_{j-1}\eta
+\alpha_{1}) .
\end{equation}
By Lemma \ref{dilation-1}, we see that
    \begin{equation}\label{proof6.1-2}
  {\rm Step}_2(\pi^{(0)})\in\overline{\mathcal{C}}_\eta(\alpha_{2},\ldots,\alpha_{\lambda};\eta,k,r|q_j,q_{j})
    \end{equation}
    and
    \begin{equation}\label{proof6.1-5}|{\rm Step}_2(\pi^{(0)})|=|\sigma|+q_{j}\eta+\alpha_1=|{\rm Step}_1(\pi^{(0)})|+q_{j}\eta+\alpha_1.\end{equation}

 Case 2: Otherwise, set
\begin{equation}\label{main-step2b}
{\rm Step}_2(\pi^{(0)})=\sigma \quad \text{and} \quad
 {\rm Step}_2(\delta^{(1)})=(q_1\eta+\alpha_1,\ldots,q_j\eta
+\alpha_1).
\end{equation}
Go to Step 3 directly.

{\bf Step 3}:   Set $\varsigma={\rm Step}_2(\pi^{(0)})$, and assume that there are  $N(\varsigma)$ parts marked with $k-1$ in $RG(\varsigma)$.
From \eqref{main-step2a} and \eqref{main-step2b}, we see that

 \[{\rm Step}_2(\delta^{(1)})=(q_1\eta+\alpha_1,\ldots,q_c\eta
+\alpha_1)
,\]
where $c=j$ if $\varsigma=\sigma$, or $c=j-1$ if $\varsigma\neq \sigma$.
We next intend to merge $q_c\eta+\alpha_1,\ldots,q_1\eta+\alpha_1$  into $\varsigma$ by  successively applying   the $(k-1)$-insertion operation.  To apply the $(k-1)$-insertion $I^{\alpha_1}_{q_{c}}$ to $\varsigma$, we are required to show that
$\varsigma$ is an overpartition in $\overline{\mathcal{B}}^{\,\alpha_1}_<(\alpha_{1},\ldots,\alpha_{\lambda};\eta,k,r|N(\varsigma),q_{c})$.
There are two cases.

Case 1: $c=j-1$. By \eqref{proof6.1-2}, we see that $\varsigma \in \overline{\mathcal{C}}_\eta(\alpha_{2},\ldots,\alpha_{\lambda};\eta,k,r|q_j,
q_j)$. Since $q_{j-1}>q_{j}$, then by   Theorem \ref{sskr1}, we deduce  that $ \varsigma \in \overline{\mathcal{B}}^{\,\alpha_1}_<(\alpha_{1},\ldots,\alpha_{\lambda};\eta,k,r|N(\varsigma),q_{j-1}).$

Case 2: $c=j$. In this case, we have $\varsigma=\sigma=\pi^{(m-j)}$ and $q_{j}\geq N(\pi^{(m-j)})=N(\varsigma)$. Let $\tilde{r}_{1}(\varsigma)>\cdots>\tilde{r}_{N(\varsigma)}(\varsigma)$ be the $(k-1)$-marked parts in $RG(\varsigma)$. Assume that $p$ is the smallest integer such that $0\leq p \leq N(\varsigma)$ and $\tilde{r}_{p+1}(\varsigma)+\eta\leq \overline{(q_j-p)\eta+\alpha_1}$. It is easy to see that such $p$ exists  since $\overline{(q_j-N(\varsigma))\eta+\alpha_1}>0\geq\tilde{r}_{N(\varsigma)+1}(\varsigma)+\eta$ with the convention that $\tilde{r}_{N(\varsigma)+1}(\varsigma)=-\infty$. Note that there are no overlined parts $\equiv \alpha_1\pmod{\eta}$ in $\varsigma$, so the largest overlined part $\equiv \alpha_1\pmod{\eta}$ in $\varsigma$ is less than $\overline{(q_j-p)\eta+\alpha_1}$.

 To show that  $\varsigma \in \overline{\mathcal{B}}^{\,\alpha_1}_<(\alpha_{1},\ldots,\alpha_{\lambda};\eta,k,r|N(\varsigma),q_{j}),$ it remains to show that  if $q_j=N(\varsigma)$ and $f_{\varsigma}(0,\eta]=r-1$, then $q_j\geq 1$ and $\tilde{r}_{N(\varsigma)}(\varsigma)\leq\eta$. Assume that $q_j=N(\varsigma)$ and $f_{\varsigma}(0,\eta]=r-1$. We consider the following two subcases.

 Subcase 2.1: $\overline{\alpha_\lambda}$ is a part of $\varsigma$. Since $q_j=N(\varsigma)$ and $\varsigma\not\in\overline{\mathcal{C}}_\lambda(\alpha_{2},\ldots,\alpha_{\lambda};\eta,k,r|N(\varsigma),N(\varsigma))$,    by the definition of $\overline{\mathcal{C}}_\lambda(\alpha_{2},\ldots,\alpha_{\lambda};\eta,k,r|N,p)$, we deduce that  $N(\varsigma)\geq 1$.
  Moreover, under the condition that  $f_{\varsigma}(0,\eta]=r-1$,  we derive that $\tilde{r}_{N(\varsigma)}(\varsigma)\leq\eta$.

 Subcase 2.2: $\overline{\alpha_\lambda}$ is not a  part of $\varsigma$. In this case, we have $j<m$. By Theorem \ref{Dilation-reduction}, we see that $f_{\pi^{(0)}}(0,\eta]=f_{\pi^{(1)}}(0,\eta]=\cdots=f_{\pi^{(m-j)}}(0,\eta]=f_{\varsigma}(0,\eta]=r-1$.
 It yields that $\overline{\alpha_\lambda}$ is a part of $\pi^{(0)}$.
 Under the assumption that $\overline{\alpha_\lambda}$ is not a  part of $\varsigma$, we find that
 $\eta$ is a non-degenerate $(k-1)$-part of $\varsigma$, and so $N(\varsigma)\geq 1$. Appealing to \eqref{proof-ss-1} in the proof of Theorem \ref{tauhh}, we obtain that $\tilde{r}_{N(\varsigma)}(\varsigma)\leq \eta$.

 Overall, we arrive at
\begin{equation*}\label{proof6.1-4}
\varsigma \in \overline{\mathcal{B}}^{\,\alpha_1}_<(\alpha_{1},\ldots,\alpha_{\lambda};\eta,k,r|N(\varsigma),q_{c}).
\end{equation*}
We proceed to merge $q_c\eta+\alpha_1,\ldots,q_1\eta+\alpha_1$  into $\varsigma$ by  successively applying   the $(k-1)$-insertion.
Denote the intermediate overpartitions  by $\varsigma^{(0)},\ldots,\varsigma^{(c)}$ with $\varsigma^{(0)}=\varsigma$ and $\varsigma^{(c)}={\rm Step}_3(\pi^{(0)})$.
 Assume that there are $N(\varsigma^{(i)})$ parts marked with $k-1$ in $RG(\varsigma^{(i)})$, where $0\leq i\leq c$.

 Set  $b=0$ and repeat the following procedure until $b=c$.

 (A) Merge $q_{c-b}\eta+\alpha_1$ into $\varsigma^{(b)}$ to generate an overlined part $\equiv \alpha_1\pmod{\eta}$.  More precisely,  applying the $(k-1)$-insertion $I^{\alpha_1}_{q_{c-b}}$ to $\varsigma^{(b)}$, we get   \[\varsigma^{(b+1)}=I^{\alpha_1}_{q_{c-b}}(\varsigma^{(b)}).\]
  By Lemma \ref{isub-lem-a}, we see that
 \begin{equation*}\label{insertbbb-n}
 \varsigma^{(b+1)} \in \overline{\mathcal{B}}^{\,\alpha_1}_=(\alpha_{1},\ldots,\alpha_{\lambda};\eta,k,r|N(\varsigma^{(b)}),q_{c-b})
 \end{equation*}
 and
 \[
 |\varsigma^{(b+1)}|=|\varsigma^{(b)}|+q_{c-b}\eta+\alpha_1.\]

 (B) Replace $b$ by $b+1$. If $b=c$, then we are done. If $b<c$, then by $q_{c-b}>q_{c-b+1}$, we conclude from Theorem \ref{ssins} that
 \[ \varsigma^{(b)} \in \overline{\mathcal{B}}^{\,\alpha_1}_<(\alpha_{1},\ldots,\alpha_{\lambda};\eta,k,r|N(\varsigma^{(b)}),q_{c-b}).\]
  Go back to (A).

The above procedure generates an overpartition ${\rm Step}_3(\pi^{(0)})=\varsigma^{(c)}$ such that
 \[{\rm Step}_3(\pi^{(0)})\in \overline{\mathcal{B}}^{\,\alpha_1}_=(\alpha_{1},\ldots,\alpha_{\lambda};\eta,k,r|N(\varsigma^{(c-1)}),q_{1})\]
     and
     \begin{equation} \label{wei-par2-n}
        |{\rm Step}_3(\pi^{(0)})|=|{\rm Step}_2(\pi^{(0)})|+(q_{c}\eta+\alpha_1)
         +\cdots+(q_1\eta+\alpha_1).
    \end{equation}
Set  $\tau={\rm Step}_3(\pi^{(0)})$.  From the construction of the $(k-1)$-insertion, it can be seen that $\tau$ is an overpartition in $\overline{\mathcal{B}}_1(\alpha_1,\ldots,\alpha_\lambda;\eta,k,r)$ with $c$ overlined parts $\equiv \alpha_1\pmod{\eta}$. Taking into account  \eqref{thm6.1aaa},  \eqref{sum-m-t}, \eqref{proof6.1-5}, along with \eqref{wei-par2-n}, we can deduce that $|\tau|=|\delta^{(1)}|+|\delta^{(\lambda)}|+|\pi|$. Therefore, $\Theta$ is the desired map from
$\mathcal{D}_{\alpha_1}\times\mathcal{D}_{\alpha_\lambda}\times \mathcal{\overline{B}}_1(\alpha_{2},\ldots,\alpha_{\lambda-1} ;\eta,
k-1,r-1)$ to $\overline{\mathcal{B}}_1(\alpha_1,\ldots,\alpha_\lambda;\eta,k,r)$.

To prove that $\Theta$ is a bijection, we shall give the description of the inverse map $\Upsilon$ of $\Theta$. Let $\tau$ be an overpartition in  $\overline{\mathcal{B}}_1(\alpha_1,\ldots,\alpha_\lambda;\eta,k,r)$.
We shall construct a triple $\Upsilon(\tau)=(\delta^{(1)},\delta^{(\lambda)}, \pi)\in\mathcal{D}_{\alpha_1}\times\mathcal{D}_{\alpha_\lambda}\times \mathcal{\overline{B}}_1(\alpha_{2},\ldots,\alpha_{\lambda-1} ;\eta,
k-1,r-1)$ such that $|\tau|=|\pi|+|\delta^{(1)}|+|\delta^{(\lambda)}|$. There are two cases.

Case 1: If there are no overlined parts $ \equiv\alpha_1\pmod\eta$ in $\tau$ and there are no non-degenerate parts of $\tau$, then set $\delta^{(1)}=\emptyset$,   set ${\delta}^{(\lambda)}$ to be the ordinary partition consisting of all the parts $\equiv \alpha_{\lambda} \pmod{\eta}$ in $\tau$, and set $\pi$ to be the overpartition  consisting of all the parts $\not\equiv\alpha_{\lambda}\pmod\eta$ in $\tau$.  Clearly, $\pi\in{\mathcal{\overline{B}}_1}(\alpha_{2},\ldots,\alpha_{\lambda-1};\eta,k-1,r-1)$ and $|\tau|=|\pi|+|\delta^{(1)}|+|\delta^{(\lambda)}|$.

Case 2: If  there are  overlined parts  $\equiv\alpha_1\pmod\eta$ in $\tau$ or there are  non-degenerate parts of $\tau$, then we  iteratively apply the $(k-1)$-separation and the $(k-1)$-subtraction to $\tau$. There are three steps. We denote  the   resulting pairs   by $( {\rm Step}_i(\delta^{(1)}), {\rm Step}_i(\tau))$ after the Step $i$, where $i=1,2,3$.

{\bf Step 1}:  Assume that there are $c\geq 0$ overlined parts $\equiv\alpha_1\pmod\eta$ in $\tau$. We  eliminate the $c$ overlined parts $\equiv\alpha_1\pmod\eta$ from $\tau$
   by successively applying the $(k-1)$-separation operation. Denote the intermediate pairs by  $(\gamma^{(0)}, \tau^{(0)}),\ldots,(\gamma^{(c)}, \tau^{(c)})$ with $(\gamma^{(0)}, \tau^{(0)})=(\emptyset, \tau)$.   There are two cases.

Case 1: $c=0$. Set $\gamma^{(c)}=\emptyset$ and $ \tau^{(c)}=\tau$.

Case 2: $c\geq 1$. Assume that $\overline{\eta t_0+\alpha_1}>\cdots>\overline{\eta t_{c-1}+\alpha_1}$ are the overlined parts $\equiv\alpha_1\pmod\eta$ in $\tau$. 

Set $b=0$ and execute the following procedure.   Assume that there are $N(\tau^{(i)})$ parts marked with $k-1$ in  $RG(\tau^{(i)})$, where $0\leq i\leq c$.

(A) Let $\hat{\tau}^{(b)}$ be the overpartition obtained from $\tau^{(b)}$ by removing $\overline{\eta t_{b}+\alpha_1}$.  Assume that there are $N(\hat{\tau}^{(b)})$ parts marked with  $k-1$
     in $G(\hat{\tau}^{(b)})$,  denoted $\tilde{g}_1(\hat{\tau}^{(b)})>\cdots>\tilde{g}_{N(\hat{\tau}^{(b)})}
    (\hat{\tau}^{(b)})$, and $p_{b}$ is the smallest integer such that $\tilde{g}_{p_{b}+1}(\hat{\tau}^{(b)})<\overline{\eta t_{b}+\alpha_1}$. Let $s^{(b)}=p_b+t_b$.
    By definition,
 \[\tau^{(b)} \in \overline{\mathcal{B}}^{\,\alpha_1}_=(\alpha_1,\ldots,\alpha_\lambda;\eta,k,r|
 N(\hat{\tau}^{(b)}),s^{(b)}).\]
Apply the $(k-1)$-separation ${J}^{\alpha_1}_{s^{(b)}}$ to $\tau^{(b)}$ to get $\tau^{(b+1)}$, that is,
 \[\tau^{(b+1)}={J}^{\alpha_1}_{s^{(b)}}(\tau^{(b)}).\]
 By means of Lemma \ref{isub-lem}, we find that $N(\hat{\tau}^{(b)})=N(\tau^{(b+1)})$,
 \[\tau^{(b+1)}\in \overline{\mathcal{B}}^{\,\alpha_1}_<(\alpha_1,\ldots,\alpha_\lambda;\eta,k,r|
 N(\tau^{(b+1)}),s^{(b)}),\]
 and
 \[  \left|\tau^{(b+1)}\right|=|\tau^{(b)}|-(s^{(b)}\eta+\alpha_1).\]
 Then insert $s^{(b)}\eta+\alpha_1$ into $\gamma^{(b)}$ to obtain $\gamma^{(b+1)}$.

(B)  Replace $b$ by $b+1$. If $b=c$, then we are done. Otherwise, go back to (A).

Set
\[{\rm Step}_1(\tau)={\tau}^{(c)} \quad \text{and} \quad {\rm Step}_1(\delta^{(1)})=\gamma^{(c)}=(s^{(0)}\eta+\alpha_1,
 \ldots,s^{(c-1)}\eta+\alpha_1).\]
 From the above construction, we see that
 \begin{equation}\label{step1-t}
 {\rm Step}_1(\tau) \in \overline{\mathcal{B}}^{\,\alpha_1}_<(\alpha_1,\ldots,\alpha_\lambda;\eta,k,r|
 N(\tau^{(c)}),s^{(c-1)})
 \end{equation}
 and
 \begin{equation}\label{weight-rel}
 |\tau|=|{\rm Step}_1(\tau)|+ |{\rm Step}_1(\delta^{(1)})|.
 \end{equation}
  Observe that for $0\leq b\leq c$,  there are $c-b$ overlined parts $\equiv\alpha_1\pmod\eta$ in $\tau^{(b)}$.
 Moreover,  there are no overlined parts $\equiv\alpha_1\pmod\eta$   in ${\rm Step}_1(\tau)=\tau^{(c)}$, so we deduce that
    \[{\rm Step}_1(\tau) \in \mathcal{\overline{C}}_1(\alpha_{2},\ldots,\alpha_{\lambda};\eta,k,r).\]
 Theorem  \ref{ssins} reveals that for $0\leq b< c-1$,
  \begin{equation*}
  s^{(b)}> s^{(b+1)}\geq N(\hat{\tau}^{(b+1)})=N({\tau}^{(b+2)}),
  \end{equation*}
that is,
 \begin{equation}\label{ttt1}
 s^{(0)}> s^{(1)}>\cdots>s^{(c-1)}\geq N({\tau}^{(c)}).
  \end{equation}

{\bf Step 2}: Set $\varsigma={\rm Step}_1(\tau)$. We consider the following   two cases.

Case 1: There is no non-degenerate $(r-1)$-part of $\varsigma$. Then set   ${\rm Step}_2(\tau)={\rm Step}_1(\tau)$ and
$ {\rm Step}_2(\delta^{(1)})={\rm Step}_1(\delta^{(1)})
$.

Case 2: There is a non-degenerate $(r-1)$-part of $\varsigma$. Let $\overline{\varsigma}$ be the  overpartition obtained  by subtracting $\alpha_1$ from the   non-degenerate $(r-1)$-part of $\varsigma$. Assume that there are $N(\overline{\varsigma})$ parts marked with $k-1$ in $G(\overline{\varsigma})$. Using the definition of $\overline{\mathcal{C}}_\eta(\alpha_{2},\ldots,\alpha_{\lambda};\eta,k,r|N,p)$, we find that
\begin{equation}\label{step2-t}
\varsigma\in\overline{\mathcal{C}}_\eta(\alpha_{2},\ldots,\alpha_{\lambda};\eta,k,r|N(\overline{\varsigma}),
N(\overline{\varsigma})).
\end{equation}
Apply the $(k-1)$-subtraction  $S_{N(\overline{\varsigma})\eta+\alpha_1}$ to $\varsigma$ to obtain ${\rm Step}_2(\tau)$. More precisely,
    \[{\rm Step}_2(\tau)=S_{N(\overline{\varsigma})\eta+\alpha_1}(\varsigma).\]
  Set
\[{\rm Step}_2(\delta^{(1)})=(s^{(0)}\eta+\alpha_1,
 \ldots,s^{(c-1)}\eta+\alpha_1,N(\overline{\varsigma})\eta+\alpha_1).\]

    It follows from Lemma \ref{2etatau} that
    \begin{equation*}\label{proofi6.1-2}
  {\rm Step}_2(\tau)\in\overline{\mathcal{C}}_\lambda(\alpha_{2},\ldots,\alpha_{\lambda};\eta,k,r|N(\overline{\varsigma}),N(\overline{\varsigma}))
    \end{equation*}
    and
\begin{equation}\label{proofi6.1-5}|{\rm Step}_2(\tau)|=|\varsigma|-(N(\overline{\varsigma})\eta+\alpha_1)=|{\rm Step}_1(\tau)|-(N(\overline{\varsigma})\eta+\alpha_1).\end{equation}
Moreover, there is no non-degenerate $(r-1)$-part of ${\rm Step}_2(\tau)$.

Combining \eqref{step1-t} and \eqref{step2-t}, and by Theorem \ref{sskr1},  we derive that if $c\geq 1$, then
\begin{equation}\label{tttt2}
s^{(c-1)}>N(\overline{\varsigma}).
\end{equation}

{\bf Step 3}: We proceed to eliminate the    non-degenerate $(k-1)$-parts of ${\rm Step}_2(\tau)$. There are two cases.

Case 1: There are no non-degenerate $(k-1)$-parts of ${\rm Step}_2(\tau)$. Then set ${\rm Step}_3(\delta^{(1)})={\rm Step}_2(\delta^{(1)})$ and ${\rm Step}_3(\tau)={\rm Step}_2(\tau)$.

Case 2: There are certain non-degenerate $(k-1)$-parts of ${\rm Step}_2(\tau)$. Denote the intermediate pairs by $(\zeta^{(0)},\sigma^{(0)}),(\zeta^{(1)},\sigma^{(1)}),$ and so on, with $(\zeta^{(0)},\sigma^{(0)})=({\rm Step}_2(\delta^{(1)}),{\rm Step}_2(\tau))$. 

Set $b=0$ and carry out the following procedure.

(A)  Let $\overline{\sigma}^{(b)}$ be the overpartition obtained by subtracting $\alpha_1$ from the smallest non-degenerate $(k-1)$-part of $\sigma^{(b)}$. Assume that  there are $N(\overline{\sigma}^{(b)})$ parts marked with $k-1$ in $G(\overline{\sigma}^{(b)})$, denoted  $\tilde{g}_1(\overline{\sigma}^{(b)})>\cdots>\tilde{g}_{N(\overline{\sigma}^{(b)})}(\overline{\sigma}^{(b)})$. Let $m_b$ be the largest integer such that $\{\sigma_{m_b+l}^{(b)}\}_{0\leq l\leq k-2}$ is a   non-degenerate $(k-1)$-band of ${\sigma}^{(b)}$.
     Set  $p^{(b)}$ to be the smallest integer such that $\tilde{g}_{p^{(b)}+1}(\overline{\sigma}^{(b)})< \sigma_{m_b}^{(b)}+\eta$. By definition, we see that
\[\sigma^{(b)}\in {\mathcal{C}}_\eta(\alpha_2,\ldots,\alpha_\lambda;\eta,k,r|N(\overline{\sigma}^{(b)}),p^{(b)}).\]
 Applying the $(k-1)$-subtraction $S_{p^{(b)}\eta+\alpha_1}$ to $\sigma^{(b)}$, we get
\[\sigma^{(b+1)}=S_{p^{(b)}\eta+\alpha_1}(\sigma^{(b)}).\]
It follows from Lemma \ref{2etatau} that
\[
\sigma^{(b+1)}\in {\mathcal{C}}_{\lambda}(\alpha_2,\ldots,\alpha_\lambda;\eta,k,r|N(\overline{\sigma}^{(b)}),p^{(b)})
\]
and
\begin{equation*}\label{aaaa}
|\sigma^{(b+1)}|=|\sigma^{(b)}|-(p^{(b)}\eta+\alpha_1).
\end{equation*}
Then insert $p^{(b)}\eta+\alpha_1$ into $\zeta^{(b)}$ to generate a new partition $\zeta^{(b+1)}$.

(B)  Replace $b$ by $b+1$. If there are no non-degenerate $(k-1)$-parts of $\zeta^{(b)}$, then we are done. Otherwise, go back to (A).

Assume that the above procedure terminates with $b=j$,  set ${\rm Step}_3(\tau)=\sigma^{(j)}$ and ${\rm Step}_3(\delta^{(1)})=\zeta^{(j)}$. More precisely,
if ${\rm Step}_1(\delta^{(1)})={\rm Step}_2(\delta^{(1)})$, then
\[{\rm Step}_3(\delta^{(1)})=(s^{(0)}\eta+\alpha_1, \ldots,s^{(c-1)}\eta+\alpha_1, p^{(0)}\eta+\alpha_1,\ldots,p^{(j-1)}\eta+\alpha_1).\]
 If ${\rm Step}_1(\delta^{(1)})\neq{\rm Step}_2(\delta^{(1)})$, then
\[{\rm Step}_3(\delta^{(1)})=(s^{(0)}\eta+\alpha_1,
 \ldots,s^{(c-1)}\eta+\alpha_1,N(\overline{\varsigma})\eta+\alpha_1, p^{(0)}\eta+\alpha_1,\ldots,p^{(j-1)}\eta+\alpha_1).\]

  By  Theorem \ref{tauhh}, we see that for $0\leq b<j-1$,
\begin{equation*}\label{ttt2}
p^{(b+1)}<p^{(b)}< N(\overline{\sigma}^{(b)}).
\end{equation*}
It implies that
\begin{equation} \label{ttt3}
N(\overline{\varsigma})\geq N(\overline{\sigma}^{(0)})>p^{(0)}>p^{(1)}>\cdots>p^{(j-1)}.
\end{equation}

Set $\delta^{(1)}={\rm Step}_3(\delta^{(1)})$. Combining \eqref{ttt1}, \eqref{tttt2} and \eqref{ttt3}, we conclude that $\delta^{(1)}$ is a partition with distinct parts $\equiv \alpha_1 \pmod{\eta}$.

From the above construction, it is easy to see that
 \begin{equation}\label{sum}
|{\rm Step}_3(\tau)|=|{\rm Step}_2(\tau)|-(p^{(0)}\eta+\alpha_1)-\cdots-(p^{(j-1)}\eta+\alpha_1).
\end{equation}
Moreover, there are no non-degenerate parts of ${\rm Step}_3(\tau)$. It means that there is a part $\equiv \alpha_\lambda \pmod{\eta}$ in each  $(k-1)$-band of ${\rm Step}_3(\tau)$. If $f_{{\rm Step}_3(\tau)}(0,\eta]=r-1$, then   $\overline{\alpha_\lambda}$ is a part of ${\rm Step}_3(\tau)$.

Let ${\delta}^{(\lambda)}$ be the ordinary partition  consisting of all the parts   $\equiv \alpha_{\lambda} \pmod{\eta}$ in ${\rm Step}_3(\tau)$ and let $\pi$ be the overpartition obtained by removing all the parts $\equiv\alpha_{\lambda}\pmod\eta$ in ${\rm Step}_3(\tau)$. It is easy to see that $\delta^{(\lambda)}\in\mathcal{D}_{\alpha_\lambda}$. Since   there are no non-degenerate parts of ${\rm Step}_3(\tau)$ and there are no parts $\equiv\alpha_{\lambda}\pmod\eta$ in $\pi$, we could deduce that  $\pi \in \overline{\mathcal{B}}_1(\alpha_{2},\ldots,\alpha_{\lambda-1};\eta,k-1,r-1)$. Combining \eqref{weight-rel}, \eqref{proofi6.1-5}  and \eqref{sum}, it is easy to check that  $|\tau|=|\pi|+|\delta^{(1)}|+|\delta^{(\lambda)}|$. Thus, we complete the proof.  \qed

\section{Example}

We provide an example for the illustration of the bijection $\Theta$ in Theorem \ref{lambdathm}.

{\noindent \bf An example for the map $\Theta$ and its inverse map $\Upsilon$:} Assume that $k=6$, $r=4$, $\lambda=4$, $\eta=10$,  $\alpha_1=1$, $\alpha_2=4$, $\alpha_3=6$ and $\alpha_4=9$. Let $\delta^{(1)}=(61,51,31,21,1)$ and $\delta^{(4)}=(49,39,29,19,9)$ and let $\pi=(50,50,50,\overline{36},\overline{34},
\overline{30},30,\overline{26},\overline{24},
{20},\overline{10},10,\overline{6})$ be an overpartition in $\overline{\mathcal{B}}_1(4,6;10,5,3)$.

We unite $\pi$ and $\delta^{(4)}$ to obtain $\pi^{(0)}$, whose reverse Gordon marking is given below.
\[\begin{split}RG(\pi^{(0)})=&({50}_1,{50}_2,{50}_3,
\overline{49}_4,\overbrace{{ \overline{39}_1, \overline{36}_2,\overline{34}_3,\overline{30}_4},{ {30}_5}}^{{ \{30\}_5}},\overline{29}_1,\overline{26}_2,\overline{24}_3,{20}_4,\\
&\  \overline{19}_1,{\overline{10}_2,{10}_3,\overline{9}_1,\overline{6}_4}).\end{split}\]
We see that $\pi^{(0)}$ is an overpartition in $\overline{\mathcal{C}}_1(4,6,9;10,6,4)$ and there are no non-degenerate parts of $\pi^{(0)}$. We wish to merge the parts of $\delta^{(1)}$
into $\pi^{(0)}$.

{\bf Step 1:} Note that $N(\pi^{(0)})=1$, so  we could apply the $5$-addition to merge $1$ into $\pi^{(0)}$. It is easy to check that $\pi^{(0)}\in\overline{\mathcal{C}}_4(4,6,9;10,6,4|1,0).$

Apply the $5$-addition $A_{1}$ to $\pi^{(0)}$ to get $\pi^{(1)}$, namely, change $\overline{39}$ to $40$.
\[\begin{split}RG(\pi^{(1)})=&(\overbrace{{ {50}_1,{50}_2,{50}_3,\overline{49}_4},{ {40}_5}}^{{ \{40\}_5}},\overbrace{{ \overline{36}_1, \overline{34}_2,\overline{30}_3,{30}_4},{ \overline{29}_5}}^{{ \{\overline{29}\}_5}},\overline{26}_1,\overline{24}_2,{20}_3,\\
&\  \overline{19}_4,{\overline{10}_1,{10}_2,
\overline{9}_3,\overline{6}_4}).\end{split}
\]

By Lemma \ref{dilation-1}, we see that $
   \pi^{(1)}\in\overline{\mathcal{C}}_{10}(4,6,9;10,6,4|
   1,0).$

Note that $N(\pi^{(1)})=2$, so the first step terminates and set ${\rm Step}_1(\pi^{(0)})=\pi^{(1)}$ and ${\rm Step}_1(\delta^{(1)})=(61,51,31,21).$

{\bf Step 2:}  Denote ${\rm Step}_1(\pi^{(0)})$ by $\sigma$.  There are two parts marked with $5$ in $RG(\sigma)$, which are $\tilde{r}_1(\sigma)=40$ and $\tilde{r}_2(\sigma)=\overline{29}$. Moreover, $f_{\sigma}(0,10]=3$,  $\tilde{r}_2(\sigma)=\overline{29}>10$ and $\overline{9}$ is a part of $\sigma$. So   $\sigma\in\overline{\mathcal{C}}_4(4,6,9;10,6,4|2,2).$

Apply the $5$-addition $A_{21}$ to $\sigma$ to get ${\rm Step}_2(\pi^{(0)})$, namely, first change  $\overline{29}$ and ${40}$ in $\sigma$ to  $\overline{39}$ and ${50}$ respectively and then change $\overline{9}$ in $\sigma$ to ${10}$. We get
\[\begin{split}RG({\rm Step}_2(\pi^{(0)}))=&({ {50}_1,{50}_2,{50}_3,{50}_4},{ \overline{49}_5},{ \overline{39}_1, \overline{36}_2,\overline{34}_3,\overline{30}_4},{ {30}_5},\overline{26}_1,\overline{24}_2,\\
&\  { {20}_3,\overline{19}_4,\overline{10}_1,{10}_2},{ {10}_5},\overline{6}_3).\end{split}\]
Set ${\rm Step}_2(\delta^{(1)})=(61,51,31)$.

{\bf Step 3:} Denote ${\rm Step}_2(\pi^{(0)})$ by $\varsigma$,  we will apply the $5$-insertion to merge $61$, $51$ and $31$ of $\delta^{(1)}$  into $\varsigma$  successively to generate some overlined parts $\equiv 1 \pmod{10}$. Let $\varsigma^{(0)}=\varsigma$.
\begin{itemize}
\item Merge $31$ into $\varsigma^{(0)}$ and set $s=3$.
 There are three parts marked with $5$ in $RG(\varsigma^{(0)})$, which are $\tilde{r}_1(\varsigma^{(0)})=\overline{49}$, $\tilde{r}_2(\varsigma^{(0)})=30$ and $\tilde{r}_3(\varsigma^{(0)})={10}$. In this occasion, $p=3$ is the smallest integer such that $\overline{(3-p)\cdot 10+1}=\overline{1}\geq \tilde{r}_{p+1}(\varsigma^{(0)})+10=-\infty$ and there are no overlined parts $\equiv 1 \pmod{10}$ in $\varsigma^{(0)}$. Hence
$\varsigma^{(0)}\in\overline{\mathcal{B}}^{\,1}_<(1,4,6,9;10,6,4|3,3).$

Apply the $5$-insertion $I^{1}_{3}$ to $\varsigma^{(0)}$ to get $\varsigma^{(1)}$. More precisely, note that $p=3$, so we first change $\overline{49}$, $30$ and $10$ to $\overline{59}$, $40$ and $20$ respectively and then insert $\overline{1}$ into the resulting overpartition.
\[\begin{split}RG(\varsigma^{(1)})=&({ \overline{59}_1,{50}_2,{50}_3,{50}_4},{ {50}_5},{ {40}_1,\overline{39}_2, \overline{36}_3,\overline{34}_4},{ \overline{30}_5}, \overline{26}_1,\overline{24}_2,\\
&\  {20}_3,{20}_4,{ \overline{19}_5},{{\overline{10}_1,{10}_2}},\overline{6}_3,\overline{1}_4).\end{split}\]
As asserted by Lemma  \ref {isub-lem-a}, we have
$\varsigma^{(1)}\in\overline{\mathcal{B}}^{\,1}_=(1,4,6,9;10,6,4|3,3).$

\item Merge $51$ into $\varsigma^{(1)}$ and set $s=5$.
There are three parts marked with $5$ in $RG(\varsigma^{(1)})$, which are $\tilde{r}_1(\varsigma^{(1)})={50}$, $\tilde{r}_2(\varsigma^{(1)})=\overline{30}$ and $\tilde{r}_3(\varsigma^{(1)})=\overline{19}$. Moreover, $p=1$ is the smallest integer such that $\overline{(5-p)\cdot 10+1}=\overline{41}\geq \tilde{r}_{p+1}(\varsigma^{(1)})+10=\overline{40}$. Given that $\varsigma^{(1)}\in\overline{\mathcal{B}}^{1}_=(1,4,6,9;10,6,4|3,3).$ Theorem \ref{ssins}  indicates that $\varsigma^{(1)}\in\overline{\mathcal{B}}^{\,1}_<(1,4,6,9;10,6,4|3,5).$

Apply the $5$-insertion $I^{1}_{5}$ to $\varsigma^{(1)}$ to get $\varsigma^{(2)}$. More precisely, note that $p=1$, so we first change  $50$  to $60$ and then insert $\overline{41}$ into the resulting overpartition.
\[\begin{split}RG(\varsigma^{(2)})=&({ {60}_1,\overline{59}_2,{50}_3,{50}_4},{ {50}_5},{ \overline{41}_1,{40}_2,\overline{39}_3, \overline{36}_4},{ \overline{34}_5},\overline{30}_1, \overline{26}_2,\overline{24}_3,\\
&\  {20}_1,{20}_4,{ \overline{19}_5},{{\overline{10}_2,{10}_3}},\overline{6}_1,\overline{1}_4).\end{split}\]
Again, by Lemma  \ref {isub-lem-a}, we have
$\varsigma^{(2)}\in\overline{\mathcal{B}}^{\,1}_=(1,4,6,9;10,6,4|3,5).$

\item  Merge $61$ into $\varsigma^{(2)}$  and set $s=6$.
There are three parts marked with $5$ in $RG(\varsigma^{(2)})$, which are $\tilde{r}_1(\varsigma^{(2)})={50}$, $\tilde{r}_2(\varsigma^{(2)})=\overline{34}$ and $\tilde{r}_3(\varsigma^{(2)})=\overline{19}$. It is easy to check that $p=0$ is the smallest integer such that $\overline{(6-p)\cdot 10+1}=\overline{61}\geq \tilde{r}_{p+1}(\varsigma^{(2)})+10=60$. Knowing that $\varsigma^{(2)}\in\overline{\mathcal{B}}^{\,1}_=(1,4,6,9;10,6,4|3,5),$ Theorem \ref{ssins} indicates that $\varsigma^{(2)}\in\overline{\mathcal{B}}^{\,1}_<(1,4,6,9;10,6,4|3,6).$

Apply the $5$-insertion $I^{1}_{6}$ to $\varsigma^{(2)}$ to get $\varsigma^{(3)}$, namely, insert $\overline{61}$ as a part into $\varsigma^{(2)}$.
\begin{equation}\label{reverse-theta}\begin{split}
RG(\varsigma^{(3)})=&(\overline{61}_1,{ {60}_2,\overline{59}_3,{50}_1,{50}_4},{ {50}_5},{ \overline{41}_2,{40}_3,\overline{39}_1, \overline{36}_4},{ \overline{34}_5},\overline{30}_2, \overline{26}_1,\overline{24}_3,\\
&\  {20}_2,{20}_4,{ \overline{19}_5},{{\overline{10}_1,{10}_3}},\overline{6}_2,\overline{1}_4).\end{split}
\end{equation}

Moreover, set ${\rm Step}_3(\pi^{(0)})=\varsigma^{(3)}$  and ${\rm Step}_3(\delta^{(1)})=\emptyset.$

\end{itemize}

Set $\tau={\rm Step}_3(\pi^{(0)})$. Clearly, $\tau$ is an overpartition in $\overline{\mathcal{B}}_1(1,4,6,9;10,6,4)$ such that $|\tau|=|\pi|+|\delta^{(1)}|+|\delta^{(4)}|$.

{\noindent \bf The inverse map $\Upsilon$:} Conversely, let $\tau$ be an overpartition in $\overline{\mathcal{B}}_1(1,4,6,9;10,6,4)$ whose reverse Gordon marking  given by  \eqref{reverse-theta}. The triple $(\delta^{(1)},\delta^{(4)},\pi)$ can be obtained by iteratively using the $5$-separation and the $5$-subtraction.

{\bf Step 1:} Note that there are three overlined parts $\equiv 1 \pmod{10}$ in $\tau$. Let $\gamma^{(0)}=\emptyset$ and $\tau^{(0)}=\tau$. We will iteratively use the $5$-separation to eliminate $\overline{61}$, $\overline{41}$ and $\overline{1}$ from $\tau$.
\begin{itemize}
\item    Eliminate $\overline{61}$ from $\tau^{(0)}$ and set $t_0=6$.

  Let $\hat{\tau}^{(0)}$ be the overpartition obtained from $\tau^{(0)}$ by removing $\overline{61}$, which has the Gordon marking
\[\begin{split}G(\hat{\tau}^{(0)})=&({ {60}_5},{ \overline{59}_1,{50}_4,{50}_3,{50}_2},\overline{41}_1,{ {40}_5},{ \overline{39}_4, \overline{36}_3,\overline{34}_2,\overline{30}_1},\overline{26}_4,\overline{24}_3,\\
&\  { {20}_5},{ {20}_2,\overline{19}_1,\overline{10}_4,{10}_3},\overline{6}_2,\overline{1}_1).\end{split}\]
There are three parts marked with $5$ in $G(\hat{\tau}^{(0)})$, which are $\tilde{g}_1(\hat{\tau}^{(0)})=60$, $\tilde{g}_2(\hat{\tau}^{(0)})={40}$ and $\tilde{g}_3(\hat{\tau}^{(0)})={20}$. Moreover,  $p_0=0$ is the smallest integer such that $\overline{t_0\cdot 10+1}=\overline{61}>\tilde{g}_{p_0+1}(\hat{\tau}^{(0)})=60$. Set $s^{(0)}=p_0+t_0=6$. Then $\tau^{(0)}\in \overline{\mathcal{B}}^{\,1}_=(1,4,6,9;10,6,4|3,6).$

We then apply the $5$-separation $J^{1}_{6}$ to $\tau^{(0)}$.  In other words, $\gamma^{(1)}$ is obtained from  $\tau^{(0)}$ by removing  $\overline{61}$, which means that $\tau^{(1)}=\hat{\tau}^{(0)}$ and $\gamma^{(1)}=(61)$. Appealing to Lemma \ref{isub-lem}, we deduce that $\tau^{(1)}\in \overline{\mathcal{B}}^{\,1}_<(1,4,6,9;10,6,4|3,6).$

\item  Eliminate  $\overline{41}$ from  $\tau^{(1)}$  and set $t_1=4$.

Let $\hat{\tau}^{(1)}$ be the overpartition obtained from $\tau^{(1)}$ by removing $\overline{41}$. We have
\[\begin{split}G(\hat{\tau}^{(1)})=&({ {60}_5},{ \overline{59}_4,{50}_3,{50}_2,{50}_1},{ {40}_5},{ \overline{39}_4, \overline{36}_3,\overline{34}_2,\overline{30}_1},\overline{26}_4,\overline{24}_3,\\
&\  { {20}_5},{ {20}_2,\overline{19}_1,\overline{10}_4,{10}_3},\overline{6}_2,\overline{1}_1).\end{split}\]
There are three parts marked with $5$ in $G(\hat{\tau}^{(1)})$, which are $\tilde{g}_1(\hat{\tau}^{(1)})=60$, $\tilde{g}_2(\hat{\tau}^{(1)})={40}$ and $\tilde{g}_3(\hat{\tau}^{(1)})=20$. Now, $p_1=1$ is the smallest integer such that $\overline{ t_1\cdot 10+1}=\overline{41}>\tilde{g}_{p_1+1}(\hat{\tau}^{(1)})=40$. Set $s^{(1)}=p_1+t_1=5$ and we get
$\tau^{(1)}\in \overline{\mathcal{B}}^{\,1}_=(1,4,6,9;10,6,4|3,5).$    Clearly,  $s^{(0)}>s^{(1)}$, in agreement with Theorem \ref{ssins}.

  Apply the $5$-separation $J^{1}_{5}$ to $\tau^{(1)}$ to get $\tau^{(2)}$. We first remove $\overline{41}$ from $\tau^{(1)}$ to get $\hat{\tau}^{(1)}$, and then change $60$ in $\hat{\tau}^{(1)}$ to  $50$ to obtain ${\tau}^{(2)}$. Finally, we insert $51$ into $\gamma^{(1)}$ to obtain $\gamma^{(2)}$. Hence  $\gamma^{(2)}=(61,51)$, and
 \[\begin{split}G({\tau}^{(2)})=&({{\overline{59}_5,{{50}_4},{50}_3,{50}_2,{50}_1}},{{{40}_5},{\overline{39}_4, \overline{36}_3,\overline{34}_2,\overline{30}_1}},\overline{26}_4,\overline{24}_3,\\
&\  {{{20}_5},{{20}_2,\overline{19}_1,\overline{10}_4,{10}_3}},\overline{6}_2,\overline{1}_1).\end{split}\]
We now have $\tau^{(2)}\in \overline{\mathcal{B}}^{\,1}_<(1,4,6,9;10,6,4|3,5)$, as expected by Lemma \ref{isub-lem}.

\item Finally, eliminate $\overline{1}$ from  $\tau^{(2)}$  and set $t_2=0$.

Let $\hat{\tau}^{(2)}$ be the  overpartition obtained from $\tau^{(2)}$ by removing $\overline{1}$, so that
\[\begin{split}G(\hat{\tau}^{(2)})=&({ \overline{59}_5},{ {50}_4,{50}_3,{50}_2,{50}_1},{ {40}_5},{ \overline{39}_4, \overline{36}_3,\overline{34}_2,\overline{30}_1},\overline{26}_3,\overline{24}_2,\\
&\  { {20}_5},{ {20}_4,\overline{19}_1,\overline{10}_3,{10}_2},\overline{6}_1).\end{split}\]
There are three parts marked with $5$ in $G(\hat{\tau}^{(2)})$, which are $\tilde{g}_1(\hat{\tau}^{(2)})=\overline{59}$, $\tilde{g}_2(\hat{\tau}^{(2)})={40}$ and $\tilde{g}_3(\hat{\tau}^{(2)})={20}$. Meanwhile, $p_2=3$ is the smallest integer such that $\overline{t_2\cdot10+1}=\overline{1}>\tilde{g}_{p_2+1}(\hat{\tau}^{(2)})=-\infty$.
Set $s^{(2)}=p_2+t_2=3$. Then $\tau^{(2)}\in \overline{\mathcal{B}}^{\,1}_=(1,4,6,9;10,6,4|3,3).$ In accordance with Theorem   \ref{ssins}, we have $s^{(1)}>s^{(2)}$.

Apply the $5$-separation $J^{1}_{3}$  to $\tau^{(2)}$. We first remove $\overline{1}$ from $\tau^{(2)}$ to get $\hat{\tau}^{(2)}$. Next, we   change $\overline{59}$, ${40}$ and ${20}$ in $\hat{\tau}^{(2)}$ to  $\overline{49}$, ${30}$ and ${10}$ respectively to obtain $\tau^{(3)}$. Finally, we insert $31$ into $\gamma^{(2)}$ to obtain $\gamma^{(3)}$. Hence  $\gamma^{(3)}=(61,51,31)$, and
\[\begin{split}G({\tau}^{(3)})=&({{{50}_5,{{50}_4},{50}_3,{50}_2,\overline{49}_1}},{{\overline{39}_5, \overline{36}_3,\overline{34}_2,\overline{30}_4,{30}_1}},\overline{26}_3,\overline{24}_2,\\
&\  {{{20}_5},{\overline{19}_1,\overline{10}_4,{10}_3,{10}_2}},\overline{6}_1).\end{split}\]
Using Lemma \ref{isub-lem}, we have $\tau^{(3)}\in \overline{\mathcal{B}}^{\, 1}_<(1,4,6,9;10,6,4|3,3).$
\end{itemize}

Set ${\rm Step}_1(\tau)={\tau}^{(3)}$ and ${\rm Step}_1(\delta^{(1)})=\gamma^{(3)}=(61,51,31).$
There are no overlined parts $\equiv 1 \pmod{10}$ in ${\tau}^{(3)}$. The fact that  ${\rm Step}_1(\delta^{(1)})$ is a partition with distinct parts reflects the claim of Theorem  \ref{ssins}.

{\bf Step 2:} Denote ${\rm Step}_1(\tau)$ by $\varsigma$. Since $f_{\varsigma}(0,10]=3$, there are no $5$-bands of $\varsigma$ in $(0,\overline{19})$  and $\overline{9}$ is not a part of $\varsigma$, we see that there is a non-degenerate $3$-part of $\varsigma$. We will apply the $5$-subtraction to $\varsigma$ to obtain an overpartition in $\overline{\mathcal{C}}_1(4,6,9;10,6,4)$ so that there is  no non-degenerate $3$-part.

 To this end,  we first subtract $1$ from a $10$ in  $\varsigma$ to obtain  $\overline{\varsigma}$. We get
\[\begin{split}G(\overline{\varsigma})=&({ {50}_5},{ {50}_4,{50}_3,{50}_2,\overline{49}_1},{ \overline{39}_5},{  \overline{36}_4,\overline{34}_3,\overline{30}_2,{30}_1},\overline{26}_4,\overline{24}_3,\\
&\  {{{20}_2},{\overline{19}_1,\overline{10}_4,{10}_3}},\overline{9}_2,\overline{6}_1).\end{split}\]
There are two $5$-marked parts in $G(\overline{\varsigma})$, which are $\tilde{g}_1(\overline{\varsigma})={50}$ and $\tilde{g}_2(\overline{\varsigma})=\overline{39}$.  Hence
$ {\varsigma}\in \overline{\mathcal{C}}_{10}(4,6,9;10,6,4|2,2).$
Next, we  change ${50}$  and $\overline{39}$ in $\overline{\varsigma}$ to ${40}$  and $\overline{29}$ respectively to obtain ${\rm Step}_2(\tau)$.
\[\begin{split}G({\rm Step}_2(\tau))=&({50}_5,{50}_4,{50}_3,\overline{49}_2,{ {40}_1}, { \overline{36}_4,\overline{34}_3,\overline{30}_2,{30}_5},\overline{29}_1,\overline{26}_4,\overline{24}_3,\\
&\  {20}_2,\overline{19}_1,\overline{10}_4,{10}_3,\overline{9}_2,\overline{6}_1).\end{split}\]
Set ${\rm Step}_2(\delta^{(1)})=(61,51,31,21)$.

{\bf Step 3:} Denote ${\rm Step}_2(\tau)$ by $\sigma$.   We will iteratively apply the $5$-subtraction  to $\sigma$ to obtain an overpartition in $\overline{\mathcal{C}}_1(4,6,9;10,6,4)$ so that there is no non-degenerate $5$-part.

The smallest non-degenerate part of $\sigma$ is $40$ in the 5-band    $\{40,\overline{36},\overline{34},\overline{30},30\}$. We first subtract $1$ from   $40$ in $\sigma$ to get $\overline{\sigma}$. We get
\begin{equation}\label{examaa}
\begin{split}G(\overline{\sigma})=&({{{50}_4,{50}_3,{50}_2,\overline{49}_1}},\overline{39}_1,\overline{36}_4,\overline{34}_3,\overline{30}_2,{ {30}_5},{  \overline{29}_1,\overline{26}_4,\overline{24}_3,{{20}_2}},\\
&\  {{\overline{19}_1,\overline{10}_4,{10}_3}},\overline{9}_2,\overline{6}_1).\end{split}
\end{equation}
There is one part marked with $5$ in $G(\overline{\sigma})$, which is $\tilde{g}_1(\overline{\sigma})={30}$. Moreover $p=0$ is the smallest integer such that $\tilde{g}_{p+1}(\overline{\sigma})={30}<50=
40+10$. Hence $\sigma\in \overline{\mathcal{C}}_{10}(4,6,9;10,6,4|1,0)$.

We then apply the $5$-subtraction $S_{1}$  to $\sigma$ to get $\sigma^{(1)}$. Then $\sigma^{(1)}$  is obtain from $\sigma$ by subtracting $1$ from   ${40}$ in $\sigma$.  Set $\zeta^{(1)}=(61,51,31,21,1)$.
 Note that there is  no non-degenerate $5$-part of $\sigma^{(1)}$, so we  set  ${\rm Step}_3(\delta^{(1)})=\zeta^{(1)}$ and ${\rm Step}_3(\tau)=\sigma^{(1)}$, where the  Gordon marking of ${\rm Step}_3(\tau)$ is given in \eqref{examaa}.

 Set
 \[\delta^{(1)}={\rm Step}_3(\delta^{(1)})=(61,51,31,21,1),\]
 and set
 \[\delta^{(4)}=(49,39,29,19,9) \quad \text{and} \quad
 \pi=(50,50,50,\overline{36},\overline{34},\overline{30},30,\overline{26},\overline{24},{20},\overline{10},10,\overline{6}),\]
 where $\delta^{(4)}$ consists of all the parts $\equiv 9 \pmod{10}$ in ${\rm Step}_3(\tau)$.

Obviously,  $\delta^{(1)}=(61,51,31,21,1)\in\mathcal{D}_1$, $\delta^{(4)}=(49,39,29,19,9)\in\mathcal{D}_9$ and $\pi$ is an overpartition in  $\overline{\mathcal{B}}_1(4,6;10,5,3)$ such that $|\tau|=|\delta^{(1)}|+|\delta^{(4)}|+|\pi|$.

 \vskip 0.2cm
\noindent{\bf Acknowledgment.} This work
was supported by  the National Natural Science Foundation of China. The first author is partially supported by
Science \& Technology Department of Sichuan Province (Nos. 24NSFSC4848).
We are grateful to the referees  for their insightful suggestions leading to an improvement of an earlier version.

\end{document}